%% file: hp3d_user_guide.tex
\documentclass[11pt]{report}

\input{Preambles/preamble_report}

\input{Preambles/preamble_stefan}

\begin{document}

\baselineskip=16pt
\parskip=4pt

%%%%%%%%%%%%%%%%%%%%%%%%%%%%%%%
% TITLE PAGE
%%%%%%%%%%%%%%%%%%%%%%%%%%%%%%%
\begin{titlepage}

\title
{
\Huge
\hp3D User Manual
}

\author
{
\LARGE
Stefan Henneking
\hskip 40pt
Leszek Demkowicz
}

\affil
{
\large \bf
Oden Institute for Computational Engineering and Sciences\\
\large \bf
The University of Texas at Austin
}

\date
{
\large
Austin, TX. \the\year{}.
}

\maketitle

\end{titlepage}

\clearpage
\pagenumbering{roman}
\setcounter{page}{2} 

%%%%%%%%%%%%%%%%%%%%%%%%%%%%%%%
% ACKNOWLEDGEMENT
%%%%%%%%%%%%%%%%%%%%%%%%%%%%%%%
\chapter*{Acknowledgements}

The authors would like to thank the developers and collaborators who have contributed to the development of the \hp3D software. In addition to the authors, the main contributors to the current version of the \hp3D finite element code are (in alphabetical order):

\begin{itemize}
	\itemsep 0pt
	\item Jacob Badger
	\item Ankit Chakraborty
	\item Federico Fuentes
	\item Paolo Gatto
	\item Brendan Keith
	\item Kyungjoo Kim
	\item Jaime D. Mora
	\item Sriram Nagaraj
	\item Socratis Petrides
\end{itemize}

Additionally, the authors would like to thank each one of the reviewers who have helped with writing and improving this user manual.

\vskip 25pt

The development of the \hp3D software and documentation, including this user manual, was partially funded by NSF award \#2103524.

%%%%%%%%%%%%%%%%%%%%%%%%%%%%%%%
% TABLE OF CONTENTS
%%%%%%%%%%%%%%%%%%%%%%%%%%%%%%%
\clearpage
\tableofcontents

%%%%%%%%%%%%%%%%%%%%%%%%%%%%%%%
% CHAPTERS
%%%%%%%%%%%%%%%%%%%%%%%%%%%%%%%
\clearpage
\pagenumbering{arabic}

% CHAPTER 1: INTRODUCTION
\graphicspath{{Figures/1_INTRODUCTION/}}
\input{Chapters/1_INTRODUCTION/introduction_chapter}

% CHAPTER 2: INSTALLING
\graphicspath{{Figures/2_INSTALLING/}}
\input{Chapters/2_INSTALLING/installing_chapter}

% CHAPTER 3: BASICS
\graphicspath{{Figures/3_BASICS/}}
\input{Chapters/3_BASICS/basics_chapter}

% CHAPTER 4: ADVANCED
\graphicspath{{Figures/4_ADVANCED/}}
\input{Chapters/4_ADVANCED/advanced_chapter}

% CHAPTER 5: PARALLEL
\graphicspath{{Figures/5_PARALLEL/}}
\input{Chapters/5_PARALLEL/parallel_chapter}

%%%%%%%%%%%%%%%%%%%%%%%%%%%%%%%
% APPENDICES
%%%%%%%%%%%%%%%%%%%%%%%%%%%%%%%
% \appendixpage
\appendix

% \titleformat{\chapter}[display]{\large \bfseries \center}{Appendix \thechapter}{20pt}{\large}
\graphicspath{{Figures/APPENDIX_A_EXAMPLES/}}
\input{Chapters/APPENDIX_A_EXAMPLES/appendix_examples_chapter}

\graphicspath{{Figures/APPENDIX_B_APPLICATIONS/}}
\input{Chapters/APPENDIX_B_APPLICATIONS/appendix_applications_chapter}

\graphicspath{{Figures/APPENDIX_C_DPG/}}
\input{Chapters/APPENDIX_C_DPG/appendix_dpg_chapter}

%%%%%%%%%%%%%%%%%%%%%%%%%%%%%%%
% REFERENCES
%%%%%%%%%%%%%%%%%%%%%%%%%%%%%%%

\phantomsection
\addcontentsline{toc}{chapter}{Bibliography}
\printbibliography

%%%%%%%%%%%%%%%%%%%%%%%%%%%%%%%
% INDEX
%%%%%%%%%%%%%%%%%%%%%%%%%%%%%%%
% \printindex

\end{document}

%% file: Preambles/preamble_report.tex
\usepackage[T1]{fontenc}
\usepackage[latin1]{inputenc}
\usepackage[english]{babel}
\usepackage[margin=1in]{geometry}
\usepackage{amsmath}
\usepackage{amssymb}
\usepackage{amsfonts}
\usepackage{authblk}
\usepackage{array}
\newcolumntype{L}[1]{>{\raggedright\let\newline\\\arraybackslash\hspace{0pt}}p{#1}}
\newcolumntype{C}[1]{>{\centering\let\newline\\\arraybackslash\hspace{0pt}}p{#1}}
\newcolumntype{R}[1]{>{\raggedleft\let\newline\\\arraybackslash\hspace{0pt}}p{#1}}
\usepackage{float}
\usepackage{graphicx}
\usepackage{caption}
\usepackage{subcaption}
\usepackage{titletoc}
\usepackage[page,toc,title,titletoc]{appendix} %page,toc,titletoc,title
\usepackage{xfrac}
\usepackage{placeins}
\usepackage{bm}
\usepackage{amsthm}
\usepackage[x11names]{xcolor}
\usepackage{comment}
\usepackage{arydshln}
\usepackage{changepage}
\usepackage{booktabs} % nice tables

\usepackage{verbatim} % verbatim environment ignores LaTeX commands

\usepackage{etoolbox}
\makeatletter
\patchcmd{\chapter}{\if@openright\cleardoublepage\else\clearpage\fi}{}{}{}
\makeatother

\usepackage[
	backend=bibtex,
	giveninits=true,
	maxbibnames=10]{biblatex}
  \definecolor{ICES}{RGB}{94, 156, 174}
  \definecolor{ORANGE}{RGB}{191, 87, 0}
  \definecolor{RED}{RGB}{190, 30, 49}
  \definecolor{SUN}{RGB}{227, 81, 51}
  \definecolor{GREEN}{RGB}{0, 171, 86}
  \definecolor{BLUE}{RGB}{11, 78, 179}
  \definecolor{BROWN}{RGB}{122, 80, 40}
  \definecolor{GREY}{RGB}{50, 50, 50}
  \definecolor{TEAL}{RGB}{0, 160, 176}
\definecolor{codegreen}{rgb}{0,0.6,0}
\definecolor{codegray}{rgb}{0.5,0.5,0.5}
\definecolor{codepurple}{rgb}{0.58,0,0.82}
\definecolor{backcolor}{rgb}{0.95,0.95,0.92}
\usepackage{listings} % used for code snippets
\usepackage{csquotes} % When using babel or polyglossia with biblatex, loading csquotes is recommended to ensure that quoted texts are typeset according to the rules of your main language.

\usepackage{colortbl}
\usepackage{tcolorbox}
\usepackage{booktabs}

\usepackage{enumitem}

\usepackage{tikz}
\usetikzlibrary{arrows.meta}

\usepackage{anyfontsize}

\usepackage[ruled]{algorithm2e}
\usepackage{algorithmic}

\usepackage{imakeidx}
\makeindex[columns=2, title=Alphabetical Index, intoc]

\usepackage[doublespacing]{setspace} % control line spacing

\usepackage{tocloft}
\makeatletter
\g@addto@macro\appendix{%
  \addtocontents{toc}{%
    \protect%
  }%
}
\g@addto@macro\appendix{%
  \addtocontents{toc}{%
    \protect%
  }%
}
\makeatother

\usepackage{titlesec}
\titleformat{\chapter}[display]{\large \bfseries \center}{Chapter \thechapter}{20pt}{\large}
\titleformat{\section}{\normalsize \bfseries \raggedright}{\thesection}{20pt}{\normalsize}
\titleformat{\subsection}{\normalsize \bfseries \raggedright}{\thesubsection}{10pt}{\normalsize}

\usepackage[hidelinks]{hyperref}
\hypersetup{colorlinks=false}

\numberwithin{equation}{chapter}
\graphicspath{{Figures/}}

\newtheorem{remark}{Remark}[chapter]

%% file: Preambles/preamble_stefan.tex
\newcommand{\var}[1]{\textcolor{SUN}{\texttt{#1}}}

\newcommand{\code}[1]{\texttt{#1}}

\newcommand{\file}[1]{\textcolor{BROWN}{\texttt{#1}}}

\newcommand{\routine}[1]{\textcolor{BLUE}{\texttt{#1}}}

\newcommand{\module}[1]{\textcolor{codegreen}{\texttt{#1}}}

\newcommand{\omp}[1]{\textcolor{TEAL}{\texttt{#1}}}

\newcommand{\tab}[1]{\ref{tab:#1}}

\def\hp3D{$hp$3D}

\newcommand{\bb}[1]{\mathbb{#1}}

\newcommand{\mf}[1]{\mathfrak{#1}}
\newcommand{\mr}[1]{\mathrm{#1}}

\newcommand{\trial}{\mathcal{U}}

\newcommand{\test}{\mathcal{V}}

\newcommand{\eps}{\varepsilon}

\newcommand{\bs}{\boldsymbol}

\newcommand{\Gammah}{\Gamma_{\hskip -1pt h}}
\newcommand{\nablah}{\nabla_{\hskip -1pt h}}

\def\grad{\nabla}

\def\curl{\nabla \times}
\def\tcurl{\text{curl}}
\def\div{\nabla \cdot}
\def\tdiv{\text{div}}

\def\hcurl{\nabla_{\hspace{-1pt}h} \times}
\def\hHcurl{H(\text{curl}, \Omega_h)}

\def\Hcurl{H(\text{curl}, \Omega)}

\def\lb{\langle}
\def\rb{\rangle}

%% file: Chapters/1_INTRODUCTION/introduction_chapter.tex
%
%!TEX root = ../../hp3D_user_guide.tex
%

\chapter{Introduction}
\label{chap:introduction}

%--------------------------------------------------------------------

The \hp3D finite element (FE) software has been developed by Prof.~Leszek Demkowicz and his students, postdocs, and collaborators at The University of Texas at Austin over the course of many years. The current version of the software is available at \url{https://github.com/Oden-EAG/hp3d}. This user manual was written to provide guidance to current and prospective users of \hp3D. We welcome any feedback about the user manual and the \hp3D code in general; please contact us at:
\begin{itemize}
	\item Prof.~Leszek Demkowicz: leszek@oden.utexas.edu
	\item Dr.~Stefan Henneking: stefan@oden.utexas.edu
\end{itemize}

\paragraph{What and who is the \hp3D code written for?}
The \hp3D code is an academic software. For this reason, \hp3D is different from other publicly available finite element codes (e.g., FEniCS, MFEM, deal.II, ...) in several important ways:

\begin{itemize}
\item Many of these modern codes are written with the aim to hide as much of the finite element machinery from the user as possible, leaving the user to only specify weak formulations, including boundary conditions, and other features of their application in abstract form. While this software design has several advantages, e.g.~fast prototyping of a new application, and may be seen as beginner-friendly, it does not at all expose the user to important FE concepts such as shape functions, Piola maps, and element-level integration. \hp3D distinguishes itself by exposing the user to these fundamental concepts and, for this reason, is conceptually aimed toward two different kinds of audiences:
\begin{enumerate}
	\item the novel FE user who is interested in learning 3D finite element coding technology; and 
	\item the advanced FE user who is interested in having lower-level access to the software making it simpler to add custom, application-specific features to the code. 
\end{enumerate}
In \hp3D, there are only minimal abstraction layers between the lower-level data structures and the upper-level user application code. On the one hand, this direct access to data structures makes it straightforward to customize FE implementations; on the other hand, we advise that only advanced users attempt modifications of data structures that are part of the library code, because there are few safeguards protecting the user from introducing bugs. When questions about library features arise, our best advice is to contact the developers before attempting modifications.

\item \hp3D includes a host of advanced finite element features such as isotropic and anisotropic $hp$-adaptivity for hybrid meshes of all element shapes, supporting conforming discretizations of the entire $H^1$--$H(\tcurl)$--$H(\tdiv)$--$L^2$ exact-sequence spaces. These features were built based on decades of experience in FE coding which went into developing optimized data structures for $hp$-adaptive finite element computation. The \hp3D data structures come with unique algorithms, including routines for constrained approximation and projection-based interpolation based on rigorous mathematical FE theory. For more information about the FE software design of the \hp3D library, we refer to \cite{hpbook2,hpbook3}. More recent additions include the support for trace variables needed for discretizations with the discontinuous Petrov--Galerkin (DPG) method, as well as support for hybrid MPI/OpenMP-parallel computation, also detailed in \cite{hpbook3}.

\item Unlike some other FE libraries, \hp3D does \emph{not} have a full-time developer team that is working to maintain or develop features for the software. Rather, the \hp3D development has depended on a small team of collaborators writing and testing newly developed features in the code. For this reason, there is no large support or developer team that can swiftly add features upon user request.

\end{itemize}

\paragraph{Is there a 2D version of the \hp3D code?} Yes, the two-dimensional $hp$2D FE software is conceptually equivalent to the 3D code with the exception that it does not support MPI-distributed parallel computation. However, we have so far not made it publicly available but we may do so in the future if there is an interest by the community. A former version of the 2D FE code was documented in:
\begin{itemize}
	\item \fullcite{hpbook}
\end{itemize}

\begin{remark}
This is a preliminary version of the user manual for the \hp3D finite element software. The user manual is still being developed and updated. In addition to this user manual, the \hp3D code is documented in the following references and references therein:

\begin{itemize}
	\item \fullcite{hpbook2}
	\item \fullcite{hpbook3}
	\item \fullcite{fuentes2015shape}
\end{itemize}

\end{remark}

%\input{Chapters/1_INTRODUCTION/comments}

%% file: Chapters/2_INSTALLING/installing_chapter.tex
%
%!TEX root = ../../hp3D_user_guide.tex
%

\chapter{Installing the \hp3D Code}
\label{chap:installing}

%--------------------------------------------------------------------

An application written for the \hp3D software is compiled in two steps. First, the \hp3D library itself is compiled; second, the particular application, which must be linked to the \hp3D library, is compiled. Changes to the application take effect by recompiling only the application, whereas changes to the underlying library source code take effect only if both the library and the application are recompiled.

Under most circumstances, the user will not need to recompile the \hp3D library since changes in the application do not affect the library source code. Therefore, most users will only compile the \hp3D library once, and from then on exclusively modify and compile the application code.

In some instances, the user may wish to change the dependencies of the library, e.g., linking to a new version of a third-party package, which will require recompilation of the library. Another situation for which the library must be recompiled is if the user chooses to toggle one of the library-wide preprocessor variables (e.g., \var{DEBUG}) or modify the compiler arguments. It is recommended that \code{`make clean'} is always executed before recompiling the \hp3D library.

The user can download the \hp3D GitHub repository via https or ssh:
\begin{itemize}
	\item via https: \code{git clone https://github.com/Oden-EAG/hp3d.git}
	\item via ssh: \code{git clone git@github.com:Oden-EAG/hp3d.git}
\end{itemize}

The main directory is then accessed by \code{`cd hp3d/trunk'}.

\begin{remark}
The \hp3D software has been deployed on a variety of UNIX-based systems. The code runs efficiently on personal laptops, including MacBooks, small workstations and clusters, all the way to large supercomputers. Currently, we are working on providing a new and more robust build system for the \hp3D library. We do not have any experience installing \hp3D on Windows-based systems and currently have no plans to add support for such systems.
\end{remark}

%%%%%%%%%%%%%%%%%%%%%%%%%%%%%%%
\section{Compiling the Library}
\label{sec:compiling}
%%%%%%%%%%%%%%%%%%%%%%%%%%%%%%%

This section provides some basic instructions for installing the code. The user should be familiar with the system in order to provide file paths to the required third-party packages/dependencies.

The \hp3D library is written exclusively in Fortran. The Fortran compiler must be compatible with the Fortran90 standard and the compiler must support the Message Passing Interface (MPI) (e.g.~OpenMPI, MPICH). The support for OpenMP threading is optional.

\subsection{Dependencies}
\label{sec:dependencies}

The configure file \file{m\_options} (see next section) must link to valid paths for external libraries. The following external libraries are used:
\begin{itemize}
	\itemsep 0pt
	\item Intel MKL [optional]
	\item X11
	\item PETSc (all following packages can be installed with PETSc)
	\item HDF5/pHDF5
	\item MUMPS
	\item Metis/ParMetis
	\item Scotch/PT-Scotch
	\item PORD
	\item Zoltan
\end{itemize}

\subsection{Configure file}
\label{sec:configure-file}

The user must create a configure file \file{trunk/m\_options}, which provides information needed by the \file{makefile}. We recommend that the user copies an existing configure file and then modifies it as needed:
\begin{enumerate}
	\itemsep 0pt
	\item{
		Copy one of the existing \file{m\_options} files from the example files in \file{hp3d/trunk/m\_options\_files/} into \file{hp3d/trunk/}. For example:\\ 
		\code{`cp m\_options\_files/m\_options\_TACC\_intel19 ./m\_options'}.
	} \item{
		Modify the \file{m\_options} file to set the correct path to the main directory:\\ 
		Set the \var{HP3D\_BASE\_PATH} to the \emph{absolute} path of the \file{hp3d/trunk/}.
	} \item{
		To compile the library, type \code{`make'} in \file{hp3d/trunk/}. Before compiling, link to the external libraries by setting the correct paths in the \file{m\_options} file---e.g.~\var{PETSC\_LIB}, \var{PETSC\_INC}---as well as setting compiler options by modifying the \file{m\_options} file as described below.
	}
\end{enumerate}

The library compilation is governed by the user-defined preprocessing flags \var{COMPLEX} and \var{DEBUG}:
\begin{itemize}
	\itemsep 0pt
	\item{ \code{\var{COMPLEX} = 0}:\\
	Stiffness matrix, load vector(s) and solution DOFs are real-valued.\\
	Code blocks within: \code{\#if C\_MODE ... \#endif}, are disabled.
	} \item{ \code{\var{COMPLEX} = 1}:\\
	Stiffness matrix, load vector(s) and solution DOFs are complex-valued.\\
	Code blocks within: \code{\#if C\_MODE ... \#endif}, are enabled.
	} \item{ \code{\var{DEBUG} = 0}:\\
	Compiler uses optimization flags and the library performs only minimal checks during the computation.\\
	Code blocks within: \code{\#if DEBUG\_MODE ... \#endif}, are disabled.
	} \item{ \code{\var{DEBUG} = 1}:\\
	Compiler uses debug flags, and the library performs additional checks during the computation.\\
	Code blocks within: \code{\#if DEBUG\_MODE ... \#endif}, are enabled.
	}
\end{itemize}
All computations in \hp3D are executed in double-precision arithmetic; that is, depending on the \var{COMPLEX} flag, variables are declared as \code{real(8)} or \code{complex(8)}, respectively. \hp3D provides a generic variable declaration \var{VTYPE}, which may be employed by the user in the application, defaulting to either real- or complex-type depending on whether the library was compiled in real or complex mode. The user should keep in mind some external library paths may also differ depending on real- and complex-valued computation. When compiling the \hp3D library, the library path is created under either \file{hp3d/complex/} or \file{hp3d/real/}, depending on the choice of preprocessing flag \var{COMPLEX}.

Support for OpenMP threading is optional and can be enabled/disabled via preprocessing flag \var{HP3D\_USE\_OPENMP}. Additional preprocessing flags for enabling/disabling third-party libraries:
\begin{itemize}
	\item{\code{\var{HP3D\_USE\_INTEL\_MKL} = 0}:\\
	Dependency on Intel MKL package is disabled.
	} \item{\code{\var{HP3D\_USE\_INTEL\_MKL} = 1}:\\
	Dependency on Intel MKL package is enabled, providing additional solver options to the user (e.g.~Intel MKL PARDISO).}
\end{itemize}

%%%%%%%%%%%%%%%%%%%%%%%%%%%%%%%
\section{Compiling an Application}
\label{sec:application}
%%%%%%%%%%%%%%%%%%%%%%%%%%%%%%%

Applications are implemented in \file{hp3d/trunk/problems}. A few applications and model problems are provided in the public \hp3D repository and can serve as examples. For example, the directory \file{problems/POISSON/GALERKIN} contains a Galerkin FE implementation for the classical variational Poisson problem. To compile and run the problem, type \code{`make'} in the application folder, i.e.,\\
\code{`cd problems/POISSON/GALERKIN; make; ./run.sh'}.\\
The implementation of this model problem is described in Appendix~\ref{sec:poisson-galerkin}. Other model problems and applications are described in Appendices~\ref{chap:examples}--\ref{chap:applications}.

%\input{Chapters/2_INSTALLING/comments}

%% file: Chapters/3_BASICS/basics_chapter.tex
%
%!TEX root = ../../hp3D_user_guide.tex
%

\chapter{Basic Features}
\label{chap:basics}

%--------------------------------------------------------------------

This chapter discusses some of the basic concepts in \hp3D which the user should be familiar with, as well as how to set up the application-specific files and routines that need to be provided by the user.

%%%%%%%%%%%%%%%%%%%%%%%%%%%%%%%
\section{Important Global Variables}
\label{sec:global-variables}
%%%%%%%%%%%%%%%%%%%%%%%%%%%%%%%

The \hp3D code defines several important global variables, some of which are system-generated and may be used within the application code, and others which are user-defined and instruct the library code.

\subsection{System-generated variables}

Several system modules (e.g.~\module{physics}) provide useful variables to the library and the user. We show some examples of these modules and their variables and encourage the user to take a look at the modules to learn about additional variables that may be of interest. Unless otherwise mentioned, the modules are located in \file{trunk/src/modules/}.

\begin{itemize}
	\item
	{
		\code{module} \module{physics}: 
		\begin{itemize}
			\item{ \var{NR\_PHYSA}:\\
			This global integer variable defines the total number of physics variables (or physics attributes) that have been requested by the user and are stored in the data structure. Note that each physics attribute may have multiple components. Moreover, a physics attribute may only be partially supported (i.e.~defined over a part of the mesh).
			} \item{ \var{NR\_COMP}\code{(i), i=1,...,\var{NR\_PHYSA}}:\\
			This global integer array defines, for each physics variable, how many components have been requested by the user. Depending on the energy space, a component may be scalar-valued ($H^1$, $L^2$) or vector-valued ($H(\tcurl)$, $H(\tdiv)$).
			} \item{ \var{NRINDEX}:\\
			This global integer defines the total number of components that have been requested by the user. That is, it is the sum of the entries in the global array \var{NR\_COMP(:)}.
		}
		\end{itemize}
	}
	\item
	{
		\code{module} \module{assembly}:
		\begin{itemize}
			\item{ \code{\var{ALOC}(i,j), \var{BLOC}(i), i,j = 1,...,\var{NR\_PHYSA}}: \\
			These two variables are provided to the user for computing local element matrices, i.e.~the user must store the element stiffness matrix in \var{ALOC} and the element load vector in \var{BLOC}. Both \var{ALOC} and \var{BLOC} are ``super-arrays'' with submatrices (blocks) defined for interactions of each physics attribute. For example, \code{\var{ALOC}(i,j)\%nrow} defines the number of rows of the $(i,j)$-th block of the stiffness matrix and \code{\var{ALOC}(i,j)\%array(:,:)} has the corresponding real-valued or complex-valued entries. The concept will become clear when looking at model problem implementations (\routine{elem} routine) with multiple physics variables. See also Section~\ref{sec:coupled-variables} on computing with coupled variables.
			}
		\end{itemize}
	} 
	\item
	{
		\code{module} \module{data\_structure3D}:
		\begin{itemize}
			\item {\var{NRELIS}, \var{NRELES}:\\
			The \module{data\_structure3D} module stores variables and arrays related to the mesh data structure (see Section~\ref{sec:data-structure}). \var{NRELIS} is the number of initial mesh elements (does not change during computation), and \var{NRELES} is the number of active mesh elements (increases with each mesh refinement).
			}
		\end{itemize}
	}
\end{itemize}

\subsection{User-defined variables}

These global variables, set directly by the user, affect the library behavior. For example, the static condensation module \module{stc}, which implements static condensation routines for local element matrices and is used by the system routine \routine{celem\_systemI} during assembly, lets the user define three variables.

\begin{itemize}
	\item
	{
		\code{module} \module{stc}:
		\begin{itemize}
			\item {\var{ISTC\_FLAG} $\in \{$\code{.true.,.false.}$\}$: \\
				Enables/disables static condensation of bubble DOFs in element matrices.}
			\item {\var{STORE\_STC} $\in \{$\code{.true.,.false.}$\}$: \\
				Enables/disables storing Schur complement factors (recomputed otherwise).}
			\item {\var{HERM\_STC} $\in \{$\code{.true.,.false.}$\}$: \\
				Enables/disables optimized linear algebra routines for Hermitian matrices.}
		\end{itemize}
	}
	\item
	{
		\code{module} \module{parameters}: 
		\begin{itemize}
			\item {\var{MAXP} $\in \{$\code{1,2, ..., 9}$\}$:\\
			Defines the maximum polynomial order of approximation allowed to be used anywhere in the mesh. This parameter is used internally for preallocating data structures and \emph{must} be changed before compiling the \hp3D library. The maximum order currently supported is \code{\var{MAXP}=9}. The user is advised to compile the library with a lower value of \var{MAXP}, depending on the maximum anticipated order used in computation, for more efficient low-order computation.
			}
			\item {\var{NRCOMS} $\in \{$\code{1,2, ...}$\}$:\\
			Number of solution component sets (i.e., copies of each variable) stored in the data structure. Multiple copies may be useful when computing time-stepping or nonlinear solutions where keeping previous solution iterates is required.\footnote{Another option for keeping extra solution components is simply defining extra variables in the physics input file. However, defining multiple copies of the entire solution component set may be more convenient for some use cases.} This parameter may be set at runtime before mesh initialization by calling the \routine{set\_parameters} routine.
			}
			\item {\var{N\_COMS}} $\in \{$\code{1,2, ..., \var{NRCOMS}}$\}$:\\
			The solution component set the user wishes to compute, assemble and solve for (if using \code{\var{NRCOMS}>1}).\footnote{Alternatively, one can use the \routine{copy\_coms} routine to transfer solution DOFs from one component set into another.} Currently, only \code{\var{N\_COMS}=1} is supported but this feature will be enabled in a future release.
			\item {\var{NRRHS}} $\in \{$\code{1,2, ...}$\}$:\\
			Number of right-hand sides (``loads'') defined for the problem. For each additional load, additional solution variables are allocated and solved for. This parameter may be set at runtime before mesh initialization by calling the \routine{set\_parameters} routine. Currently, only \code{\var{NRRHS}=1} is supported but this feature will be enabled in a future release.
		\end{itemize}
	}
	\item
	{
		\code{module} \module{paraview}: 
		\begin{itemize}
			\item {\var{VLEVEL} $\in \{$\code{0,1,2,3,4}$\}$:\\
			Determines how refined the mesh output is when exporting geometry and solution data to ParaView/vtk. If \code{\var{VLEVEL}=0}, the mesh is outputted ``as is'' with the geometry and solution linearly interpolated for the ParaView plot. Each ``upscaling'' of \var{VLEVEL} increases the resolution by one refinement level. See Section~\ref{sec:VTK} for further details and additional user-defined ParaView module variables.
			}
		\end{itemize}
	}
	\item
	{
		\code{module} \module{physics}:
		\begin{itemize}
			\item{ \code{\var{PHYSAd}(i)} $\in \{$\code{.true.,.false.}$\},$ \code{i=1,...,\var{NR\_PHYSA}}:\\
			Used for indicating a \emph{homogeneous} Dirichlet BC for any components of the variable defined as a Dirichlet component. For such components, the computation is then more efficient because no Dirichlet data needs to be interpolated.
			}
			\item{ \code{\var{PHYSAi}(i)} $\in \{$\code{.true.,.false.}$\},$ \code{i=1,...,\var{NR\_PHYSA}}:\\ 
			Used for computing with traces (see Section~\ref{sec:traces}).
			}
			\item{ \code{\var{PHYSAm}(i)} $\in \{$\code{.true.,.false.}$\},$ \code{i=1,...,\var{NR\_PHYSA}}:\\ 
			Used for enabling/disabling physics attributes (see Section~\ref{sec:coupled-variables}).
			}
		\end{itemize}
	}
\end{itemize}

%%%%%%%%%%%%%%%%%%%%%%%%%%%%%%%
\section{Application-Specific Input Files and Routines}
\label{sec:simple-example}
%%%%%%%%%%%%%%%%%%%%%%%%%%%%%%%

\hp3D's library is set up to read certain user-defined inputs when initializing the program (e.g.~initial geometry file) and call certain user-defined routines (e.g.~providing element-local matrices). This section briefly reviews the input files and routines that the user must provide.

In general, when beginning a new application implementation, it is advisable to start off by copying an existing model implementation and then modify the required routines. Applications should be coded within the subdirectory \file{trunk/problems/}. Examples of specific model problem implementations are given in Appendix~\ref{chap:examples}.

\subsection{Input files}

Each of the three required input files---\file{control}, \file{physics}, and \file{geometry}---must be provided with a certain formatting and define a number of required parameters. For additional details and examples of the input formatting, we refer to the model problems.

\begin{itemize}
	\item{ 
		\file{control}\\
		Sets global variables in \code{module} \module{control}. The user must provide the (relative) path to this file to the \var{FILE\_CONTROL} variable in \code{module} \module{environment} before mesh initialization.
	}\item{
		\file{physics}\\
		Sets global variables in \code{module} \module{physics}. The user must provide the (relative) path to this file to the \var{FILE\_PHYS} variable in \code{module} \module{environment} before mesh initialization.
	}\item{
		\file{geometry}\\
		Defines the initial geometry mesh. The user must provide the (relative) path to this file to the \var{FILE\_GEOM} variable in \code{module} \module{environment} before mesh initialization.
	}
\end{itemize}

\subsection{Driver and required routines}

\begin{itemize}
	\item{
		\code{program} \routine{main}: \\
		The \routine{main} program is the ``driver'' of the application. It should set environment variables, initialize the mesh data structure, and interact with the user (or execute a pre-defined job). The initialization of an application includes calling various library routines, such as \routine{read\_control}, \routine{read\_input}, \routine{read\_geometry}, and \routine{hp3gen}.
	}
	\item{
		\code{subroutine} \routine{set\_initial\_mesh(\var{Nelem\_order})}: \\
		The initialization routine \routine{hp3gen} sets up the mesh data structures (see next section). During this initialization, a user-provided routine \routine{set\_initial\_mesh} is called which defines the supported physics variables and initial order of approximation for each element (stored in \var{Nelem\_order(:)}), as well as the boundary condition flags for physics components on element faces.
	}
	\item{
		\code{subroutine} \routine{dirichlet}: \\
		If the user has specified non-homogeneous Dirichlet flags for any component on the boundary, then the user-provided \routine{dirichlet} routine is required. The routine is called from the system routine \routine{update\_Ddof} which computes the Dirichlet data via projection-based interpolation. \routine{dirichlet} takes a physical coordinate input $x$ and must return the value and first derivatives of the expected $H^1$, $H(\tcurl)$, and/or $H(\tdiv)$ boundary data at coordinate $x$.
	}
	\item{
		\code{subroutine} \routine{elem}: \\
		The \routine{elem} routine is at the heart of the application code---it implements the variational formulation on the element level. This routine is called by system routines during the finite element assembly and provides the element-local stiffness matrix and load vector for each element in the active mesh.
	}
\end{itemize}

%%%%%%%%%%%%%%%%%%%%%%%%%%%%%%%
\section{Mesh Data Structure}
\label{sec:data-structure}

Due to the support for hybrid meshes and adaptive refinements, the \hp3D mesh data structure is dynamically build as refinements are executed. Before the main data structure arrays (\var{ELEMS} and \var{NODES}) are discussed, we review how element data is stored and accessed for elements of different shapes.

\subsection{Element data}
\label{sec:element-data}

\paragraph{Module \module{element\_data}.}
The \hp3D code supports four types of elements: hexahedra, tetrahedra, prisms, and pyramids. Any element computations are executed in an object-oriented programming fashion using general element utilities provided in \file{src/modules/}\module{element\_data}. For instance, a simple loop through element vertices is executed by

% Use lstlisting environment for code snippets
\begin{lstlisting}[caption=Loop over element vertices., label={lst:loop_element_vertices}]
do iv=1,NVERT(elem_type)
   ...
enddo
\end{lstlisting}

\noindent where \var{elem\_type} is a member variable of an element specifying the element type---\code{BRIC}, \code{TETR}, \code{PRIS}, \code{PYRA}---and \var{NVERT(:)} is a parameterized array (defined in \code{module} \module{node\_types}) returning the number of vertices for each element type. In this way, we avoid writing four separate versions of the code for the four types of elements but cover all four cases with just one piece of code. The user may want to review the \code{module} \module{element\_data} to learn about the existing element utilities and the underlying logic.

The \module{element\_data} module starts by listing master element coordinates for the four element vertices. This defines the geometry of master elements and establishes enumeration of master element vertices. The element edges are listed next by providing numbers of the corresponding endpoint vertices. This again serves a double purpose: we enumerate the element edges and provide the corresponding local edge orientations. Finally, we list element faces by providing face-to-vertex node connectivities. This again implies the local face orientations. Finally, following the specified orientations, we provide local parametrizations for element edges and faces. All these utilities are necessary for computing element matrices. The \routine{elem} routines provided in the model problem implementations (Appendix~\ref{chap:examples}) illustrate how to use the utilities provided in the module.

\paragraph{Initial mesh generation. Data structure arrays \var{ELEMS} and \var{NODES}.}

Following the input of geometry for the Geometry Modeling Package (GMP), an initial mesh is generated. The initial mesh generator represents an interface between the GMP package and the \hp3D code, and it does essentially two things: it generates initial mesh vertex, edge, face and element middle nodes in the data structure array \var{NODES}, and it generates initial mesh elements in data structure array \var{ELEMS}, including element-to-nodes connectivities. Both \var{ELEMS} and \var{NODES} are arrays of objects provided by the \module{data\_structure3D} module. Array \var{ELEMS} is static; its dimension equals the number of initial mesh elements equal to the number of GMP blocks. Array \var{NODES} is dynamic; it grows during the mesh refinements as new nodes are generated. The elements arising from mesh refinements are logically identified with their middle nodes. Middle node numbers are thus unique identifiers for an element, and the number of a node always coincides with the location in the \var{NODES} array. For the initial mesh, the number of an initial mesh element in \var{ELEMS} coincides with the corresponding middle node number in \var{NODES}.

\subsection{Looping through active mesh elements}
\label{sec:element-nodes}

The \var{NODES} array stores all (abstract) element nodes---vertices, edges, faces, and middle nodes---whether they are part of the current mesh or not. In order to access active mesh elements (associated with active middle nodes), the \hp3D code provides a functionality to loop through the active mesh elements in their ``natural'' order defined by the tree structure of the \var{NODES} array.

\paragraph{Natural order of elements.}

A typical loop through the active elements in the mesh is executed by using routine \file{datstrs/}\routine{nelcon}:

\begin{lstlisting}[caption=Loop over active mesh elements., label={lst:loop_active_mesh_elements}]
mdle = 0
do iel=1,NRELES
	call nelcon(mdle, mdle)
	...
enddo
\end{lstlisting}

where \var{NRELES} is a global variable storing the number of active mesh elements. The \routine{nelcon} routine is based on nodal trees restricted to element middle nodes only. See \cite{hpbook,hpbook2} for additional information.

\begin{remark}
For MPI/OpenMP parallel execution of loops over active mesh elements, the \hp3D \module{data\_structure3D} module builds two arrays, \var{ELEM\_ORDER} and \var{ELEM\_SUBD}, based on \routine{nelcon}. \var{ELEM\_ORDER} is an array of all active mesh elements, and \var{ELEM\_SUBD} is an array of active mesh elements within a subdomain of a distributed mesh. Each one is built (resp.~updated) by using \routine{nelcon} after each mesh refinement or mesh repartitioning. See Section~\ref{sec:mpi-modules} for additional details.
\end{remark}

\subsection{Uniform mesh refinements}
\label{subsec:uniform-refinements}

An isotropic uniform $h$-refinement can be executed via routine \routine{global\_href}. Anisotropic uniform $h$-refinement routines are also provided for some cases (e.g.~\routine{global\_href\_aniso\_bric}) but should be used carefully since particular anisotropic refinements may not be compatible with a hybrid mesh of different element types (see refinement flags in Section~\ref{sec:adaptive-refinements}).

Analogously, isotropic uniform $p$-refinements are easily executed via routine \routine{global\_pref}. For the case of $p$-refinements, the user may also execute global unrefinements (i.e.~uniformly lowering the polynomial order of approximation) via routine \routine{global\_punref}.

Adaptive $hp$-refinements are discussed in Section~\ref{sec:adaptive-refinements}.

%%%%%%%%%%%%%%%%%%%%%%%%%%%%%%%
\section{Fundamental Finite Element Algorithms}
\label{sec:FE_algorithms}
%%%%%%%%%%%%%%%%%%%%%%%%%%%%%%%
The code supports two fundamental FE algorithms: solution of a boundary-value problem and adaptive solution of a boundary-value problem. In the second case, the user must provide an additional routine providing an a-posteriori error estimate. This section outlines the assembly and linear solve of a finite element system in \hp3D. Adaptive solutions are discussed in Section~\ref{sec:adaptivity} (adaptive solver and refinements) and Section~\ref{sec:parallel-adaptivity} (parallel computation of an adaptive solution).

\subsection{Assembly}

The finite element assembly process takes element-local stiffness matrices and load vectors (provided by the user routine \routine{elem}) and assembles the local blocks into a sparse global matrix. This process also includes much of the sophisticated \hp3D machinery for constrained approximation, as well as modifications due to accounting for Dirichlet data and static condensation of element-interior degrees of freedom. The assembly routines are hidden from the user and do not usually need to be interacted with directly. Instead, the user calls one of the provided linear solver interfaces which perform the assembly process for the user application.

The global stiffness matrix and load vector are assembled automatically by the solver interface, i.e.~the assembly routines loop through all {\em active} elements in the mesh and compute the corresponding element matrices. This implementation of the element stiffness matrix and load vector are in the user-provided element routine \routine{elem}. For each active element \var{mdle}, the code calls \routine{elem} from routine \file{src/constrs/}\routine{celem\_system} which returns the information about the {\em modified element} corresponding to \var{mdle} and the modified element matrices. For a detailed discussion of \routine{celem\_system}, see \cite{hpbook,hpbook2}. In summary, the \routine{celem\_system} routine compiles the information about the modified element nodes, the corresponding number of DOFs, performs partial assembly to account for the constrained nodes, and enforces the Dirichlet BCs by computing the modified load vector and eliminating the Dirichlet DOFs from the system of equations. The solution DOFs returned by the solver are then stored in the data structure arrays.

\subsection{Linear solver}

The code provides interfaces to several linear solvers:
\begin{itemize}
	\item{ MUMPS: \\
		Different MUMPS solver options are available to the user via the following solver interfaces:
		\begin{itemize}
			\item{ OpenMP MUMPS: \routine{mumps\_sc} \\
			Sequential MUMPS solver with OpenMP support.
			}
			\item{ MPI/OpenMP MUMPS: \routine{par\_mumps\_sc} \\
			Distributed MUMPS solver with OpenMP support.
			}
			\item{ MPI/OpenMP Nested Dissection: \routine{par\_nested} \\
			Statically condenses subdomain interior DOFs onto subdomain interfaces and subsequently solves the coupled interface problem via the distributed MUMPS solver.
			}
		\end{itemize}
	}
	\item{ MKL\_PARDISO: \routine{pardiso\_sc} \\
		Intel MKL's PARDISO solver (OpenMP support only) is available via solver interface \routine{pardiso\_sc}, provided the library code was compiled with flag \code{\var{HP3D\_USE\_INTEL\_MKL} = 1} (see Section~\ref{sec:configure-file}).
	}
	\item{ Frontal Solver: \routine{solve1} \\
		Homegrown sequential solver (no MPI or OpenMP support) used mainly for debugging purposes. The solver interface is implemented in the routine \routine{solve1}.
	}
	\item{ PETSc: \routine{petsc\_solve} \\
		The PETSc solver interface supports a variety of external solver options, depending upon the PETSc library \hp3D was linked to during compilation (cf.~Section~\ref{sec:dependencies}).
	}
\end{itemize}
Some of the solver interfaces take a single-character input argument \code{`H'} or \code{`G'}, indicating whether the linear system is symmetric/Hermitian or not, in order to perform optimized linear algebra computations. In some cases, the user may also specify \code{`P'} to indicate positive-definiteness of the system.

The auxiliary data structures for assembly and linear solve of the finite element system are deallocated before the solver interface returns. On return, the finite element solution (i.e.~solution degrees of freedom) has been stored in the data structure and is available to the user for post-processing.

%%%%%%%%%%%%%%%%%%%%%%%%%%%%%%%
\section{Exporting Geometry Mesh and Solution to VTK}
\label{sec:VTK}
%%%%%%%%%%%%%%%%%%%%%%%%%%%%%%%

\hp3D supports two output formats for exporting geometry meshes and solutions to VTK: XDMF and VTU. Each one has certain compatibilities with output options that are indicated in Table~\tab{paraview-compatibility}. The geometry and solution data are stored as HDF5 (XDMF) or raw binary files (VTU). The user interfaces for exporting mesh and solution data are provided in the \module{paraview} module and several related subroutines. Note that exporting meshes with pyramids is currently not supported but will be supported in a future release of the code.

\paragraph{ParaView module.}
\hp3D's internal ParaView module, \module{src/modules/paraview}, defines the data structures and parameters required for writing ParaView/VTK output data. To the user, most of the functionality of this module stays hidden, except for various parameters that may be defined by the user. These user-defined parameters include:
\begin{itemize}
	\item \var{PARAVIEW\_DIR}\\
	Path to output directory; can be an absolute path or a path relative to the binary execution file. The output path must be valid (the module does not create new directories).
	\item \var{VIS\_VTU} $\in \{$ \code{.true.}, \code{.false.} $\}$\\
	Enables/disables VTU output format. By default, the output format is set to XDMF.
	\item \var{PARAVIEW\_DUMP\_GEOM} $\in \{$ \code{.true.}, \code{.false.} $\}$\\
	Controls whether a geometry file should be written by the ParaView output routines. For example, writing a geometry file on every call to the ParaView driver may be useful if successive solutions are visualized for an adaptive mesh. On the other hand, if visualizing a transient solution on a fixed mesh, the geometry file only needs to be written once at the initial time step. Note that if using the VTU format, this parameter is not supported and the mesh data is always written.
	\item \var{PARAVIEW\_DUMP\_ATTR} $\in \{$ \code{.true.}, \code{.false.} $\}$\\
	Controls whether solution data for physics attributes should be written by the ParaView output routines. In some circumstances, the user may want to output only the geometry data to visualize the finite element mesh.
	\item \var{VLEVEL} $\in \{$ \code{`0', `1', `2', `3', `4'} $\}$\\
	By default, \hp3D's utilities output ParaView files for visualizing linear elements; in this case, higher output resolution in \hp3D is achieved by artificially ``upscaling'' the element data. This upscaling is done by evaluating geometry and solution data at additional points within each element, corresponding to visualizing a uniformly $h$-refined element. The level of upscaling is determined by the value of the \module{paraview} module parameter \var{VLEVEL}. At most four levels of upscaling are supported.
	\item \var{SECOND\_ORDER\_VIS} $\in \{$ \code{.true.}, \code{.false.} $\}$\\
	Enables/disables writing second-order elements (both geometry and solution data). Note that by default, ParaView uses linear interpolation between the Lagrange points. To enable high-order interpolation in ParaView, the user needs to increase the ``Nonlinear Subdivision Level.'' If \var{SECOND\_ORDER\_VIS} is enabled, then additional upscaling is not supported (i.e., \code{\var{VLEVEL}=`0'}).
	\item \var{PARAVIEW\_TIME}\\
	If set by the user, this real-valued parameter is used as a timestamp for the current solution output.
	\item \code{\var{PARAVIEW\_LOAD}(i)} $\in \{$\code{.true.,.false.}$\},$ \code{i=1,...,\var{NRRHS}}\\
	Enables/disables writing a specific solution vector. The user is advised to set the values by calling \routine{paraview\_select\_load}. By default, all entries are enabled.
	\item \code{\var{PARAVIEW\_ATTR}(i)} $\in \{$\code{.true.,.false.}$\},$ \code{i=1,...,\var{NR\_PHYSA}}\\
	Enables/disables writing a specific physics variable. The user is advised to set the values by calling \routine{paraview\_select\_attr}. By default, all entries are enabled.
	\item \code{\var{PARAVIEW\_COMP\_REAL}(i)} $\in \{$\code{.true.,.false.}$\},$ \code{i=1,...,\var{NRINDEX}}\\
	Enables/disables writing the real part of a specific component. The user is advised to set the values by calling \routine{paraview\_select\_comp\_real}. By default, all entries are enabled.
	\item \code{\var{PARAVIEW\_COMP\_IMAG}(i)} $\in \{$\code{.true.,.false.}$\},$ \code{i=1,...,\var{NRINDEX}}\\
	Enables/disables writing the imaginary part of a specific component. The user is advised to set the values by calling \routine{paraview\_select\_comp\_imag}. By default, all entries are enabled.
\end{itemize}

\begin{table}[htb]
\centering
\caption{Compatibility of ParaView options with output formats}
\label{tab:paraview-compatibility}
\begin{tabular}{lll}
\toprule
Variable & XDMF format & VTU format \\
\midrule
\var{PARAVIEW\_DUMP\_GEOM} & yes & no (always writes mesh data) \\
\var{VLEVEL} & yes & yes \\
\var{SECOND\_ORDER\_VIS} & no mixed meshes & yes \\
\bottomrule
\end{tabular}
\end{table}

\paragraph{ParaView driver.}
The ParaView driver routine serves as an interface between the user application and the \hp3D ParaView module. The main purpose of the driver routine is to coordinate the writing of output files, i.e., writing an XML or VTU metadata file that defines the ParaView fields, initiating the writing of the geometry data, and iterating through physics variables calling various ParaView routines to assemble and write the corresponding solution data. If using the VTU output format, a ParaView PVD file is written as a metadata file with timestamp information for each VTU file. With XDMF, timestamp values are written directly into the XML metadata files.

The default ParaView driver can be found in \file{src/graphics/paraview/paraview\_driver.F90}. This \routine{paraview\_driver} routine can be called directly by the user and will, by default, write geometry and solution data as defined by the variables specified above. This driver routine may also serve as a user-customizable template for more application-specific needs.

%\input{Chapters/3_BASICS/comments}

%% file: Chapters/4_ADVANCED/advanced_chapter.tex
%
%!TEX root = ../../hp3D_user_guide.tex
%

\chapter{Advanced Topics}
\label{chap:advanced}

%--------------------------------------------------------------------

This is a preliminary version of the user manual. This chapter on the advanced FE concepts provided by the \hp3D code is currently under development and will be further expanded and completed in a future version of the user manual.

%%%%%%%%%%%%%%%%%%%%%%%%%%%%%%%
\section{Custom Boundary Conditions}
\label{sec:advanced-BC}
%%%%%%%%%%%%%%%%%%%%%%%%%%%%%%%

In the preceding chapter, it was shown how to set up basic Dirichlet boundary conditions for a variable. In this section, some custom options for setting up various types of boundary conditions are discussed. Note that the \hp3D code allows the user to customize boundary conditions in an almost arbitrary way. A few common boundary condition types are presented with the goal of introducing the general idea of customizing boundary conditions for an application.

In the user-defined \routine{set\_initial\_mesh} routine, boundary conditions are set \emph{per component:}
\begin{itemize}
	\item{
	\var{ELEMS(iel)\%bcond(1:6,1:NRINDEX)} $\in \{$ \code{0,1, ..., 9} $\}$:\\
	user-defined value indicating for faces of each initial mesh element \var{iel} whether a component is a free component (\code{0}), a Dirichlet component (\code{1}), or has a custom BC (\code{2, ..., 9}).
	} \item{
	\var{NODES(nod)\%bcond(1:NRINDEX)} $\in \{$ \code{0,1} $\}$:\\
	system-generated flag indicating for each node \var{nod} whether a component is a Dirichlet component (\code{1}) or not (\code{0}).
	}
\end{itemize}
The \hp3D code provides a routine to simplify setting boundary conditions flags for particular components on (parts of) the mesh boundary:
\begin{itemize}
	\item{
		\routine{set\_bcond}: This routine sets a user-defined BC flag for a particular physics attribute component on all exterior faces of the mesh that match a certain boundary ID number. Boundary IDs for faces are defined in the mesh input file and generated when reading the mesh. If no boundary IDs are specified, then by default all boundary IDs are \code{0}.
	}
\end{itemize}

\subsection{Neumann boundary condition}
\label{sec:impedance}

Discussion of Neumann BC will be included in a future version of the user manual.

\subsection{Impedance boundary condition}
\label{sec:impedance}

Discussion of impedance BC will be included in a future version of the user manual.

%%%%%%%%%%%%%%%%%%%%%%%%%%%%%%%
\section{Adaptivity}
\label{sec:adaptivity}
%%%%%%%%%%%%%%%%%%%%%%%%%%%%%%%

The sophisticated $hp$-adaptive capabilities are one of the main selling points of the \hp3D finite element code. This section describes the basics of an adaptive solution procedure and how the user can execute $hp$-adaptive mesh refinements in the application code.

For ultraweak DPG implementations, an effective automatic $hp$-adaptation capability is developed and described in \cite{chakraborty2023hp}. Further details on using this \hp3D adaptation module will be provided in a future version of the manual.

\subsection{Adaptive solver}
\label{sec:adaptive-solver}

The simplest adaptive algorithm executes the standard logical sequence:
\[
	\text{Solve} \longrightarrow
	\text{Estimate} \longrightarrow
	\text{Mark} \longrightarrow
	\text{Refine} 
\]
until a prescribed global error tolerance is met. Error estimation involves a loop through elements and computation of element {\em error indicators} plus possible additional information aiming at selection of an optimal element refinement. Once the element error indicators are known for all active elements, selected elements are {\em marked} for pre-selected $hp$-refinements. Two marking strategies are most popular: {\em the greedy strategy} and the {\em D\"orfler strategy}. 

In the greedy strategy, a maximum element error indicator \var{error\_max} is first determined. All elements which exceed a certain percentage of \var{error\_max} are then marked for refinement. The refinement criterion for a single element is thus
\begin{center}
	\var{error} $>$ \var{perc} * \var{error\_max}.
\end{center}

In D\"orfler's strategy \cite{dorfler1996marking}, the elements are first organized in the order of descending element error indicators and the total error \var{error\_glob} is computed by summing up the element error indicators. The first $n$ elements whose cumulative sum exceeds a certain percentage of the total error are then marked for refinement, i.e.:
\[
	\sum_{j=1}^n \var{elem\_error}_j > \var{perc} * \var{error\_glob} .
\]
Marked elements are placed on a separate list along with the requested refinement flags.

\subsection{Adaptive refinements}
\label{sec:adaptive-refinements}

This section gives an overview of the $hp$-adaptive capabilities of the \hp3D code. Note that the section does \emph{not} discuss the implementation of an error indicator function for guiding mesh adaptivity, but rather introduces the mechanism of executing the adaptive refinements. It is assumed that some error indicator function, usually problem-dependent, is provided by the user. 
% Maybe mention the DPG error indicator and examples where it's used

The starting point for adaptive refiements is an array \var{elem\_ref(:)} of size \var{nr\_elem\_ref} storing a list of the middle node numbers (\var{mdle}) corresponding to elements marked for refinement.

\paragraph{$h$-adaptivity.}
In the first example, we show how $h$-adaptive refinements are executed. In essence, each middle node is refined one-by-one by calling the \routine{refine} routine with a corresponding refinement flag \var{kref} encoding the type of $h$-refinement. For middle nodes that can be anisotropically refined, \var{kref} encodes how to refine the corresponding element. Encoding of refinement flags depends upon the element type; for instance,
\begin{itemize}
	\item a hexahedral element may be refined in 7 different ways; for example:
	\begin{itemize}
		\itemsep 0pt
		\item \code{\var{kref} = 111} : refine in $x,y,z$
		\item \code{\var{kref} = 110} : refine in $x,y$
		\item \code{\var{kref} = 101} : refine in $x,z$
	\end{itemize}
	\item a prismatic element may be refined in 3 different ways:
	\begin{itemize}
		\itemsep 0pt
		\item \code{\var{kref} = 11} : refine in $xy,z$
		\item \code{\var{kref} = 10} : refine in $xy$
		\item \code{\var{kref} = 01} : refine in $z$
	\end{itemize}
\end{itemize}

For hybrid meshes, isotropic refinements can be executed by obtaining isotropic refinement flags for any element type from the routine \routine{get\_isoref} for the respective middle node. After executing $h$-adaptive refinements, it is essential that the user calls the \routine{close\_mesh} routine to preserve 1-irregularity of the mesh.

\begin{lstlisting}[caption=Isotropic $h$-adaptive refinements., label={lst:adaptive_refinement_element_loop4}]
!..iterate over elements marked for refinement
   do iel=1,nr_elem_ref
      mdle = elem_ref(iel)
!  ...get isotropic h-refinement flag for element type
      call get_isoref(mdle, kref)
!  ...h-refine element
      call refine(mdle,kref)
   enddo
!..enforce 1-irregular mesh
   call close_mesh
!..update geometry and Dirichlet DOFs
   call update_gdof
   call update_Ddof
\end{lstlisting}

\paragraph{$p$-adaptivity.}
In the second example, we show how $p$-adaptive refinements are executed. Isotropic $p$-adaptive refinements are carried out easily by the routine \routine{execute\_pref} which must only be provided with the list of elements to be refined in the form of the corresponding middle nodes. After executing adaptive $p$-refinements, the user may want to call one of the system routines \routine{enforce\_min\_rule} or \routine{enforce\_max\_rule} to apply the minimum or maximum rule, respectively, to handle neighboring elements of different approximation order in a particular way (i.e.~either lowering or raising the order of approximation on the corresponding element interfaces).

\begin{lstlisting}[caption=Isotropic $p$-adaptive refinements., label={lst:adaptive_refinement_element_loop4}]
!..execute isotropic p-refinements for a list of elements
   call execute_pref(elem_ref,nr_elem_ref)
!..enforce minimum rule
   call enforce_min_rule
!..update geometry and Dirichlet DOFs
   call update_gdof
   call update_Ddof
\end{lstlisting}

Similar to the $h$-adaptive case where the refinement flag encoded possibly anisotropic refinements, the polynomial order $p$ of element nodes encodes the possibly anisotropic order of approximation. In order to execute anisotropic $p$-adaptive refinements, the user may instead call the routine \routine{perform\_pref}, which additionally must be provided with a list of $p$-refinement flags for middle nodes corresponding to the list of elements to be refined. For each middle node on this list, the $p$-refinement flag encodes the desired order of approximation for the middle node. For example, $p=222$ encodes isotropic quadratic order for the middle node of a hexahedral element, whereas $p=34$ encodes anisotropic mixed order for the middle node of a prism. The routine will then execute corresponding $p$-refinements of edge and face nodes. After finishing $p$-adaptive refinements, the user may again want to enforce the minimum or maximum rule.

%%%%%%%%%%%%%%%%%%%%%%%%%%%%%%%
\section{Trace Variables}
\label{sec:traces}
%%%%%%%%%%%%%%%%%%%%%%%%%%%%%%%

This section discusses how trace variables can be discretized in \hp3D. This feature is especially important for discretizations with the DPG method \cite{demkowicz2017dpg}. The \hp3D code offers conforming discretizations for traces of the exact-sequence energy spaces, i.e.~functions belonging to the $H^{1/2}$-, $H^{-1/2}(\tcurl)$-, and $H^{-1/2}$-energy spaces. In \hp3D, these trace unknowns---which are defined on element interfaces---are discretized by using restrictions of $H^1$-, $H(\tcurl)$-, and $H(\tdiv)$-conforming elements to the element boundary. The discretized trace unknowns are thus continuous, tangentially continuous, and normally continuous, respectively. For further reading, we refer to \cite{demkowicz2018spaces,demkowicz2020fem}.

To realize trace variables as restrictions of exact-sequence-conforming elements, the user specifies for each physics attribute whether it is a trace variable via a global flag:

\code{\var{PHYSAi}(i)} $\in \{$\code{.true.,.false.}$\}$, \code{i=1,...,\var{NR\_PHYSA}}.

\noindent
If enabled, this implies that interactions of bubble (interior) DOFs of an element are not stored in the element matrices for the corresponding traces. In order to compute with trace unknowns, the user \emph{must} enable \hp3D's static condensation module \module{stc} (see Section~\ref{sec:global-variables}). By default, all physics variables are defined as standard variables unless otherwise instructed by the user. The user \emph{may not} for obvious reasons request a trace of an $L^2$-conforming (discontinuous) variable.

%%%%%%%%%%%%%%%%%%%%%%%%%%%%%%%
\section{Coupled Variables}
\label{sec:coupled-variables}
%%%%%%%%%%%%%%%%%%%%%%%%%%%%%%%

In \hp3D, problems with coupled variables may be defined in a way where the coupling happens over the entire domain, over part of the domain, or at an interface between different parts of the domain. This section primarily focuses on problems with variables supported in the entire mesh.

\paragraph{Computing with multiple physics variables.}

Solving multiphysics problems in \hp3D requires an understanding of how \hp3D internally stores and computes interactions between multiple variables.
First, recall the following essential information from previous sections:
\begin{itemize}
	\item Global variables \var{NR\_PHYSA}, \var{NR\_COMP(:)}, and \var{NRINDEX} store the total number of physics attributes, the number of components per physics attribute, and the total number of components, respectively.
	\item In the \routine{set\_initial\_mesh} routine, boundary conditions are set separately for each component.
	\item In the \routine{elem} routine, element-local matrices are assembled in blocks corresponding to (supported) physics attributes:\\ 
	\var{ALOC(1:NR\_PHYSA,1:NR\_PHYSA)\%array(:,:)} (element stiffness); and \\
	\var{BLOC(1:NR\_PHYSA)\%array(:,:)} (element load).\\
	Note that the load ``vector'' may have multiple columns if solving a problem for multiple right-hand sides.
\end{itemize}

Physics variables are always defined and stored (in \module{physics}) in the order of the exact sequence. For each approximation space, the user may define multiple physics attributes (in the \file{physics} input file), and each attribute may have one or multiple components (also defined in the \file{physics} input file).

The \routine{elem} routine that assembles each of the stiffness and load blocks must provide the system assembly routines with information about which blocks have been locally computed. This is necessary for two reasons:
\begin{enumerate}
	\item An element may not support all physics variables thus will only assemble a subset of the blocks corresponding to the supported attributes.
	\item A linear system for only a subset of the supported variables is assembled and solved.
\end{enumerate}
In order to provide this information to the system routines, \routine{elem} returns a flag for each of the blocks indicating whether it was computed by the element. The following arguments are defined in the \routine{elem} routine:
\begin{itemize}
	\itemsep 0pt
	\item{ \code{integer, intent(in) :: \var{Mdle}} \\
		\var{Mdle} is the unique middle node number associated with the element.
	}
	\item{ \code{integer, intent(out) :: \var{Itest}(\var{NR\_PHYSA})} \\
		For the $i$-th test function (corresponding to the $i$-th physics attribute), \code{\var{Itest}(i) = 1} indicates that the $i$-th block (row) of \var{ALOC} and \var{BLOC} corresponds to an active and supported variable.
	}
	\item{ \code{integer, intent(out) :: \var{Itrial}(\var{NR\_PHYSA})} \\
		For the $i$-th trial function (corresponding to the $j$-th physics attribute), \code{\var{Itrial}(j) = 1} indicates that the $j$-th block (column) of \var{ALOC} corresponds to an active and supported variable.
	}
\end{itemize}

The $(i,j)$-th block of \var{ALOC} is provided by the \routine{elem} routine if and only if \code{\var{Itest}(i) = 1} \emph{and} \code{\var{Itrial}(j) = 1}, and the $i$-th block of \var{BLOC} is provided if and only if \code{\var{Itest}(i) = 1}. Usually, when the $(i,j)$-th stiffness block is computed, then the $(j,i)$-th block is also provided.

\paragraph{Solving for a subset of the physics variables.}

In a problem where the user would like to compute a subset of the variables at a time (e.g.~staggered solve in simple fixed-point iteration), this can easily be done by toggling the global flags \var{PHYSAm(:)} of the \module{physics} module. For the $i$-th physics attribute (\code{i = 1,\ldots,\var{NR\_PHYSA}}), \code{\var{PHYSAm}(i)=.false.} deactivates the corresponding physics attribute. If all physics variables are otherwise supported on the entire mesh, then the \routine{elem} routine can very easily set the assembly flags \var{Itest} and \var{Itrial} as follows:

\begin{lstlisting}[caption=Solving for a subset of variables in \routine{elem}.]
!  in:   Mdle     - element middle node number
!  out:  Itest    - index for assembly
!        Itrial   - index for assembly
subroutine elem(Mdle, Itest,Itrial)
   integer, intent(in)    :: Mdle
   integer, intent(out) :: Itest(NR_PHYSA)
   integer, intent(out) :: Itrial(NR_PHYSA)
!  .......
   Itest(1:NR_PHYSA) = 0; Itrial(1:NR_PHYSA) = 0
!
!..all variables are supported on all elements,
!  thus assemble for all of the currently enabled variables
   do i = 1,NR_PHYSA
      if (PHYSAm(i)) then
         Itest(i) = 1; Itrial(i) = 1
      endif
   enddo
!  ...... compute the corresponding blocks of ALOC and BLOC
\end{lstlisting}

\begin{remark}
More information about partial support of variables (defined in parts of the domain) will be included in a future version of the user manual.
\end{remark}

%\input{Chapters/4_ADVANCED/comments}

%% file: Chapters/5_PARALLEL/parallel_chapter.tex
%
%!TEX root = ../../hp3D_user_guide.tex
%

\chapter{Parallel Computation}
\label{chap:parallel}

%--------------------------------------------------------------------

This chapter gives a brief introduction to parallel computation with MPI and OpenMP in \hp3D. Further details about parallel algorithms and data structures in \hp3D are given in \cite{hpbook3}.

%%%%%%%%%%%%%%%%%%%%%%%%%%%%%%%
\section{Modules and Variables for Parallel Computation}
\label{sec:mpi-modules}
%%%%%%%%%%%%%%%%%%%%%%%%%%%%%%%

This section briefly reviews essential functionality and variables provided in some of \hp3D's parallelization modules that are used internally and can also be leveraged by the user in the problem implementation.

\paragraph{\module{mpi\_wrapper} module.}
The \module{mpi\_wrapper} module provides two essential routines: \routine{mpi\_w\_init} and \routine{mpi\_w\_finalize}. In every problem implementation, these two routines should be the very first and the very last call, respectively, of the program:
\begin{lstlisting}[caption=Initiating and finalizing \hp3D's MPI environment., label={lst:mpi_w_init}]
program main
   use mpi_wrapper 
   call mpi_w_init
!  ... program execution ...
   call mpi_w_finalize
end program main
\end{lstlisting}
These calls are required for initializing and closing the MPI environment, including hybrid MPI/OpenMP threading, as well as initializing and closing the Zoltan environment used for load balancing (see Section~\ref{sec:load-balancing}). The calls to \routine{mpi\_w\_init} and \routine{mpi\_w\_finalize} should be made even if the program is executed by a single MPI process. 

The user is also encouraged to use the \routine{mpi\_w\_handle\_err(\var{Ierr}, \var{Str})} routine of the \module{mpi\_wrapper} module which prints the error code and a user-defined error string \var{Str} if the return code \var{Ierr} of an MPI function is not equal to \var{MPI\_SUCCESS}. For example:
\begin{lstlisting}[caption=Checking and printing MPI error codes., label={lst:mpi_w_handle_error}]
	call MPI_BARRIER (MPI_COMM_WORLD, ierr)
	call mpi_w_handle_err (ierr, "MPI_BARRIER returned with error code.")
\end{lstlisting}

\paragraph{\module{mpi\_param} module.} This module stores global parameters and variables related to MPI parallelism. The user may want to make use of the following parameters that are initialized with the MPI environment:
\begin{itemize}
	\item \code{integer, save      ::~\var{NUM\_PROCS}} \\
	Total number of MPI processes. Return value of \routine{MPI\_COMM\_SIZE} for the \var{MPI\_COMM\_WORLD} communicator.
	\item \code{integer, save      ::~\var{RANK}} \\
	MPI rank of each MPI process. \var{RANK} $\in \{ \code{0}, \code{1}, \dots, \code{\var{NUM\_PROCS}-1} \}$.
	Return value of \routine{MPI\_COMM\_RANK} for the \var{MPI\_COMM\_WORLD} communicator.
	\item \code{integer, parameter ::~\var{ROOT}  = 0} \\
	\hp3D defines the MPI process with \code{\var{RANK} 0} as the \emph{root process} (or equivalently, \emph{host process}). Frequently, the root process takes a different execution path from the remaining processes (e.g.~in user-interactive computation---see Section~\ref{sec:master-worker-paradigm}, or in I/O operations).
\end{itemize}

\paragraph{\module{par\_mesh} module.}
Most of the essential mesh (re-)partitioning functionality in \hp3D is provided by the module \file{src/modules/}\module{par\_mesh}. By default, the initial mesh is first generated on each process, then distributed according to the natural element order. The user can change the initial mesh partitioner via \routine{zoltan\_w\_set\_lb} in the \routine{initialize} routine before calling \routine{hp3gen} where the initial mesh data structure is generated and distributed. The \module{par\_mesh} module provides a \routine{distr\_mesh} routine that can be called by the user to (re-)distribute the mesh at any point, e.g., after mesh refinements. That is, the \routine{distr\_mesh} routine serves as a routine for repartitioning (resp.\ load balancing)---see Section~\ref{sec:load-balancing}.

The \module{par\_mesh} module defines two important global variables that keep track of the current state of the mesh in the parallel environment:
\begin{itemize}
	\item \code{logical, save ::~\var{DISTRIBUTED} \vspace{-5pt}}
	\begin{itemize} \itemsep 2pt
		\item \code{.false.}~if the current mesh is not distributed (all DOFs present on all MPI processes);\\
		\code{.true.~}~if the current mesh is distributed (DOFs are distributed across MPI processes).
		\item If more than one MPI process is used, the \var{DISTRIBUTED} flag should always be \code{.true.}\ after mesh initialization (once \routine{hp3gen} returns). Once distributed, a mesh cannot be ``undistributed.'' In other words, even if the partitioning is completely unbalanced (e.g., all DOFs are on one MPI process), the mesh is still distributed in the sense that each mesh element is associated with one particular subdomain (resp.\ MPI process).
		\item Initially, \code{\var{DISTRIBUTED} = .false.}\\
		After calling \routine{distr\_mesh}, \code{\var{DISTRIBUTED} = .true.}; by default, \routine{distr\_mesh} is called in routine \routine{hp3gen} during mesh generation when \code{\var{NUM\_PROCS} > 1}. 
		\item If using sequential computation (single MPI process), the value of \var{DISTRIBUTED} is irrelevant.
	\end{itemize}
	\item \code{logical, save ::~\var{HOST\_MESH} \vspace{-5pt}}
	\begin{itemize} \itemsep 2pt
	\item \code{.true.}~if and only if the \code{ROOT} process has the entire mesh (all DOFs are present on \code{ROOT}).
	\item Initially, \code{\var{HOST\_MESH} = .true.}\ since all DOFs present on all processes including the \code{ROOT} process. After calling \routine{distr\_mesh}, \code{\var{HOST\_MESH} = .false.} unless \code{\var{NUM\_PROCS} = 1}.
	\item The DOFs of a distributed mesh can be collected on the host process by calling the \routine{collect\_dofs} routine in \file{src/mpi/par\_aux.F90}. Collecting DOFs on host enables utilities written only for sequential computing (e.g., interactive visualization) and can serve as a debugging tool. After calling \routine{collect\_dofs}, \code{\var{HOST\_MESH} = .true.}
	\item If using sequential computation (single MPI process), the value of \var{HOST\_MESH} is always \code{.true.}
	\end{itemize}
\end{itemize}

As will be discussed in more detail in Section~\ref{sec:leveraging}, the user should be familiar with the meaning of these two variables to leverage parallel computation in \hp3D applications. In particular, using these mesh state variables greatly simplifies writing and executing an \hp3D application that supports both the sequential and distributed-memory setting.

\paragraph{\var{ELEM\_ORDER} and \var{ELEM\_SUBD}.}
The \module{data\_structure3D} module provides two arrays that are frequently used in the parallel environment. Recall from Section~\ref{sec:data-structure} that iterating over the active mesh elements is commonly done by using the \routine{nelcon} routine that provides a natural ordering of elements (see Listing~\ref{lst:loop_active_mesh_elements}). In the parallel setting, this sequential looping is replaced with a loop over one of two active mesh element arrays: \var{ELEM\_ORDER} and \var{ELEM\_SUBD}. The \var{ELEM\_ORDER} array stores all active mesh elements in their natural order; the \var{ELEM\_SUBD} array stores all active mesh elements within a particular subdomain of a distributed mesh. After each mesh refinement or mesh repartitioning, the two arrays \var{ELEM\_ORDER} and \var{ELEM\_SUBD} are automatically updated by the \routine{update\_ELEM\_ORDER} routine of the \module{par\_mesh} module:

\begin{lstlisting}[caption=Updating the data structure arrays \var{ELEM\_ORDER} and \var{ELEM\_SUBD}., label={lst:update_elem_order}]
mdle = 0; NRELES_SUBD = 0
do iel=1,NRELES
   call nelcon(mdle, mdle)
   ELEM_ORDER(iel) = mdle
   if (NODES(mdle)%subd .eq. RANK) then
      NRELES_SUBD = NRELES_SUBD + 1
      ELEM_SUBD(NRELES_SUBD) = mdle
   endif
enddo
if (.not. DISTRIBUTED) then
   ELEM_SUBD(1:NRELES) = ELEM_ORDER(1:NRELES)
   NRELES_SUBD = NRELES
endif
\end{lstlisting}

The total number of active mesh elements (size of \var{ELEM\_ORDER}) is denoted by \var{NRELES}, and the number of elements within the subdomain (size of \var{ELEM\_SUBD}) is denoted by \var{NRELES\_SUBD}. The \var{ELEM\_ORDER} and \var{ELEM\_SUBD} arrays are used internally for parallelizing the element loops, and they can be leveraged by the user for the same purpose as demonstrated in Section~\ref{sec:generic-element-loops}.

%%%%%%%%%%%%%%%%%%%%%%%%%%%%%%%
\section{Leveraging Parallelism in Applications}
\label{sec:leveraging}
%%%%%%%%%%%%%%%%%%%%%%%%%%%%%%%

This section discusses how MPI and OpenMP parallelism can be leveraged in the user application. As previously discussed in this chapter, most of the functionality that enables parallel computation in \hp3D is hidden from the user. Nonetheless, the user should be to some extent familiar with the parallel environment to be able to implement or modify customizable parallel algorithms within the application.

Before reading this section, it is highly recommended to review the essential modules and data structures enabling parallelism introduced in the preceding section (Section~\ref{sec:mpi-modules}).

%--------------------------------------------------------------------------------------------%
\subsection{Generic element loops}
\label{sec:generic-element-loops}

Looping over mesh elements is one of the most common operations in a finite element code. Usually, these parallelized loops are hidden from the user in \hp3D; for example, the assembly process only requires the user to provide local element matrices which are then modified and assembled to a global linear system by \hp3D's parallelized system routines. However, for any meaningful FE computation, it is essential for the user to understand how a parallel element loop works. For instance, post-processing solution data commonly requires volumetric or boundary integrals over each element. In this type of computation, most of the work is element-local thus can be effectively parallelized. Usually, the parallelized element-local computation is followed by a communication step (e.g., a reduction operator) to calculate aggregate values of interest (e.g., a global residual).

\paragraph{Subdomain element loops.}
In the distributed setting, the first level of parallelism---MPI distributed-memory computation---is implied by the mesh partitioning. Section~\ref{sec:mpi-modules} introduced the data structures and concepts needed to write a parallel loop over mesh elements within subdomains of the distributed mesh. A typical subdomain element loop can be implemented as follows:

\begin{lstlisting}[caption=Looping over active mesh elements in a subdomain., label={lst:subdomain_element_loop}]
   integer :: iel, mdle
   real(8) :: elem_val, subd_val, total_val
!..if computing sequentially (single MPI process), subdomain is the entire mesh
   if (.not. DISTRIBUTED) ELEM_SUBD(1:NRELES) = ELEM_ORDER(1:NRELES)
!..iterate over elements within the subdomain
   subd_val = 0.d0
   do iel=1,NRELES_SUBD
      mdle = ELEM_SUBD(iel)
!  ...element-local computation
      elem_val = ...
!  ...aggregate element values over subdomain
      subd_val = subd_val + elem_val
   enddo
!..aggregate subdomain values globally
   call MPI_REDUCE(subd_val,total_val,1,MPI_REAL8, &
                     MPI_SUM,ROOT,MPI_COMM_WORLD, ierr)
\end{lstlisting}

\paragraph{Parallelizing subdomain element loops.}
With MPI and OpenMP enabled, the user will typically write a loop over elements within each subdomain using OpenMP threading for parallel element computation within the subdomain. This second level of parallelism---OpenMP shared-memory computation---must be initiated explicitly by opening an OpenMP parallel environment. In this case, element-local variables must be declared as \omp{PRIVATE} variables by the user (i.e., each thread has its own copy of the variable), and subdomain aggregate values are computed by an OpenMP \omp{REDUCTION} clause. In other words, with two levels of parallelization, aggregate values are computed by a two-level reduction operation: first, a reduction across OpenMP threads within the subdomain computation; and second, a reduction across MPI processes for the global aggregate value.

\begin{lstlisting}[caption=Accelerating subdomain element loops by OpenMP threading., label={lst:omp_subdomain_element_loop}]
!$OMP PARALLEL DO PRIVATE(mdle,elem_val)    &
!$OMP REDUCTION(+:subd_val)
   do iel=1,NRELES_SUBD
      mdle = ELEM_SUBD(iel)
!  ...element-local computation
      elem_val = ...
!  ...aggregate element values over subdomain
      subd_val = subd_val + elem_val
   enddo
!$OMP END PARALLEL DO
\end{lstlisting}

\paragraph{OpenMP thread scheduling.}
By default, an OpenMP-parallel loop uses \emph{static thread scheduling}, i.e.~each loop iteration (resp.~element workload) is statically assigned a-priori to a particular OpenMP thread.
In context of $hp$-adaptively refined meshes with elements of different shape, this default strategy can lead to a severely unbalanced workload between OpenMP threads due to highly variable element workload. In this case, it is preferable to use \emph{dynamic thread scheduling}: by passing the OpenMP clause \omp{SCHEDULE(DYNAMIC)}, each thread is assigned one loop iteration (resp.~element workload) at a time. Load balancing is thus also done on two levels: the distributed-memory level (a-priori load balancing by mesh repartitioning) and the shared-memory level (dynamic load balancing by OpenMP thread scheduling). While dynamic thread scheduling at runtime eliminates the need for a-priori load balancing work, it comes with a scheduling overhead at loop execution.
This overhead can be somewhat reduced by assigning a few elements at a time (i.e., using an OpenMP \omp{chunk} size larger than one in the dynamic schedule) or by using a \omp{guided} schedule with variable chunk size. These options reduce scheduling overhead by somewhat compromising on load balance. The optimal choice of OpenMP schedule will depend on a variety of problem-dependent factors. Generally speaking: if the workload is similar for all elements, use either a static schedule or dynamic schedule with large chunk sizes; if the element workload varies much, use a dynamic schedule with small chunk sizes.

\begin{remark}
As previously mentioned, element loops for matrix assembly and other tasks are part of the system routines and thus hidden from the user application. In \hp3D, these system routines by default employ a dynamic OpenMP thread schedule in all element loops to account for varying element workload in adaptive solutions.
\end{remark}

%--------------------------------------------------------------------------------------------%
\subsection{Adaptive solution with a distributed mesh}
\label{sec:parallel-adaptivity}

A typical use-case of performing element-loops is a-posteriori error estimation. An element-local error indicator can be computed independently and in parallel for each element. For instance, in the DPG method the element-local residual serves as a built-in error indicator \cite{demkowicz2017dpg}. The following steps illustrate how parallel MPI/OpenMP computation is used for computing adaptive refinements based on D\"orfler's marking strategy
(see Section~\ref{sec:adaptive-solver}) using element residuals on a distributed mesh in \hp3D:
\begin{enumerate}
\item
	Compute element residuals within the subdomain:
	\begin{lstlisting}[caption=Adaptive refinements based on element residuals: (a) subdomain computation., label={lst:adaptive_refinement_element_loop1}]
   subd_res = 0.d0; elem_res(1:NRELES) = 0.d0
!$OMP PARALLEL DO PRIVATE(mdle,subd)    &
!$OMP SCHEDULE(DYNAMIC) REDUCTION(+:subd_res)
   do iel=1,NRELES
      mdle = ELEM_ORDER(iel)
      call get_subd(mdle, subd)
      if (RANK .ne. subd) cycle
      call elem_residual(mdle, elem_res(iel))
      subd_res = subd_res + elem_res(iel)
   enddo
!$OMP END PARALLEL DO
\end{lstlisting}
\item
	Communicate residual values by MPI reduction:
	\begin{lstlisting}[caption=Adaptive refinements based on element residuals: (b) global communication., label={lst:adaptive_refinement_element_loop2}]
   call MPI_ALLREDUCE(subd_res,total_res,1,MPI_REAL8, &
                        MPI_SUM,MPI_COMM_WORLD, ierr)
   call MPI_ALLREDUCE(MPI_IN_PLACE,elem_res,NRELES,MPI_REAL8, &
                        MPI_SUM,MPI_COMM_WORLD, ierr)
\end{lstlisting}
\item
	Mark elements for adaptive refinement (D\"orfler's strategy):
	\begin{lstlisting}[caption=Adaptive refinements based on element residuals: (c) element marking., label={lst:adaptive_refinement_element_loop3}]
!..sort elements by residual values (in descending order)
   mdle_ref(1:NRELES) = ELEM_ORDER(1:NRELES)
   call qsort_duplet(mdle_ref,elem_res,NRELES)
!..mark elements for refinement based on marking strategy
   alpha = 0.5d0 ! (marking coefficient)
   nr_elem_ref = 0; sum_res = 0.d0
   do iel=1,NRELES
      nr_elem_ref = nr_elem_ref + 1
      sum_res = sum_res + elem_res(iel)
      if (sum_res > alpha*total_res) exit
   enddo
   \end{lstlisting}
\item
	Refine marked elements:
	\begin{lstlisting}[caption=Adaptive refinements based on element residuals: (d) element refinement., label={lst:adaptive_refinement_element_loop4}]
!..iterate over elements marked for refinement
   do iel=1,nr_elem_ref
      mdle = mdle_ref(iel)
!  ...get isotropic h-refinement flag for element type
      call get_isoref(mdle, kref)
!  ...h-refine element
      call refine(mdle,kref)
   enddo
!..enforce 1-irregular mesh
   call close_mesh
!..update geometry and Dirichlet DOFs
   call update_gdof
   call update_Ddof
!..repartition the mesh for load balancing
   call distr_mesh
	\end{lstlisting}
\end{enumerate}

\FloatBarrier
%--------------------------------------------------------------------------------------------%
\subsection{Master--worker paradigm}
\label{sec:master-worker-paradigm}

There are numerous of ways of realizing a parallel driver for \hp3D applications. The user is, of course, free to design the driver of a parallel implementation as needed for the particular application of interest. This section introduces one concept that has been employed successfully within many \hp3D problem implementations. Generally, we distinguish between two different modes of operation (resp.\ program execution):

\begin{itemize}
	\item Interactive mode: \\
	A user-interactive job passes only the required arguments for initiating the program (primarily the geometry and physics input files) during the program call. Mesh refinements, adaptive solution process, visualization and I/O are all controlled interactively by the user input. This mode is suitable for sequential (single MPI process) computation or parallel computation with a moderate number of MPI processes executed on a single compute node. Usually, the interactive mode is the primary mode of execution during code development and debugging.
	\item Production mode: \\
	A production job passes all of the program arguments needed for execution during the program call. Pre-processing, mesh refinements, adaptive solutions, and post-processing operations are all driven by a pre-defined job script. This mode is suitable for any number of MPI processes and MPI/OpenMP configuration, including large-scale execution on many compute nodes.
\end{itemize}

As an example, we refer to the driver file \file{main.F90} of the Galerkin Poisson implementation in \file{problems/POISSON/GALERKIN}. In this driver implementation, a program argument ``\code{-job \var{JOB}}'' controls whether the program is executed in the interactive mode (\code{\var{JOB}.eq.0}) or production mode (\code{\var{JOB}.ne.0}). In production mode, the value of \var{JOB} may then define which one of a number of pre-defined job scripts should be executed. 

The remainder of this section is concerned with the implementation of the so-called \emph{master--worker paradigm} for the interactive mode of \hp3D applications in the parallel distributed-memory setting.

\paragraph{Interactive mode in the sequential setting.}
First, we introduce a typical interactive I/O screen employed by \hp3D's subroutines (e.g., see \file{problems/POISSON/GALERKIN/main.F90}):

\begin{lstlisting}[caption=Interactive mode in shared-memory execution., label={lst:interactive_mode_sequential}]
!..user interface displaying menu in a loop
   idec = 1
   do while(idec /= 0)
!  ...print user options
      call print_menu
!  ...read user input
      read(*,*) idec
!  ...act on user input
      select case(idec)
         case(0) ; exit
!     ...call a routine handling the user input
         call default ; exec_case(idec)
      end select
   enddo
\end{lstlisting}

For context, a basic user menu could include the following user options:

\begin{lstlisting}[caption=A basic interactive-mode user menu., label={lst:interactive_mode_menu}]
 =-=-=-=-=-=-=-=-=-=-=-=-=-=-=-=-=-=-=-=-=
 QUIT....................................0
 
         ---- Visualization ----        
 HP3D graphics (graphb)..................1
 HP3D graphics (graphg)..................2
 Paraview................................3
 
         ---- DataStructure ----       
 Display arrays (interactive)...........10
 Display DataStructure info.............11
 
          ---- Refinements ----            
 Single uniform h-refinement............20
 Single uniform p-refinement............21
 Refine a single element................22
 
            ---- Solvers ----              
 MUMPS..................................30
 Pardiso................................31
 Frontal................................32
 PETSc..................................33
 
        ---- Error / Residual ----        
 Compute exact error....................40
 Compute residual.......................41
 =-=-=-=-=-=-=-=-=-=-=-=-=-=-=-=-=-=-=-=-=
 \end{lstlisting}

Similar I/O screens are utilized in various system debugging routines the user has access to (e.g., displaying nodal information via the interactive routine \routine{src/datstrs/result}). The mechanism for executing these routines in the distributed-memory setting is similar to the one described in the context of the application driver \routine{main} discussed here.

\paragraph{Interactive mode in the parallel setting.}
As a general rule, the user always interacts with one distinct MPI process at a time. This MPI process---called the \emph{master}---is the \var{ROOT} process. All other MPI ranks are referred to as \emph{worker} processes. In this master--worker paradigm, user inputs are handled by the master process who communicates the input to the worker processes. For that reason, the parallel driver that handles the user interaction is split into two driver routines: \routine{master\_main} and \routine{worker\_main}.
\begin{lstlisting}[caption=Splitting master and worker execution paths., label={lst:master_worker_main}]
   if (JOB .ne. 0) then
!  ...execute production mode
      call exec_job
   else
!  ...execute interactive mode
      select case(RANK)
         case(ROOT) ; call master_main
         case default ; call worker_main
      end select
   endif
\end{lstlisting}

Similar to the sequential case, the \routine{master\_main}---executed only by the master process---displays an interactive user menu and waits for the user input. Once the user input is read, the master broadcasts the user input to the workers. Consequently, the \routine{worker\_main} routine---executed by all worker processes---primarily consists of a receiving broadcast (from \var{ROOT}) followed by executing the user command:

\begin{lstlisting}[caption=Interactive mode in distributed-memory execution., label={lst:interactive_mode_parallel}]
  idec=1
  do while(idec /= 0)
!  ...waiting for the master to broadcast user input
      call MPI_BCAST(idec,1,MPI_INTEGER,ROOT,MPI_COMM_WORLD, ierr)
!  ...act on user input
      select case(idec)
         case(0) ; exit
!     ...call a routine handling the user input
         call default ; exec_case(idec)
      end select
   enddo
\end{lstlisting}

In essence, this summarizes how interactive user inputs are typically handled in parallel \hp3D applications. The master--worker paradigm can, in principle, be executed for an arbitrary number of processes. Primarily, though, the interactive mode serves the user to check intermediate results of the computation interactively during code development and debugging.

\begin{remark}
In some circumstances, the user may want to interact directly with one of the worker processes. For instance, when the user wants to display information that the master process does not have access to, such as geometry degrees of freedom outside of its subdomain. In this case, the master can prompt a user input for choosing a particular worker rank to interact with rather than relaying all user commands to the worker by MPI communication.
\end{remark}

\begin{lstlisting}[caption=Initiating an interactive worker routine by a master broadcast., label={lst:interactive_mode_worker}]
!..ask the user which MPI process should display information
   if (RANK .eq. ROOT)
      write(*,*) 'Select processor RANK: '; read(*,*) r
   endif
   call MPI_BCAST(r,1,MPI_INTEGER,ROOT,MPI_COMM_WORLD, ierr)
!..open interactive data structure routine from worker process with rank r
   if (r .eq. RANK) call result
!..wait for the interactive worker routine to finish
   call MPI_BARRIER(MPI_COMM_WORLD, ierr)
\end{lstlisting}

%%%%%%%%%%%%%%%%%%%%%%%%%%%%%%%
\section{Load Balancing}
\label{sec:load-balancing}
%%%%%%%%%%%%%%%%%%%%%%%%%%%%%%%

The basic algorithm for an adaptive solution with dynamic load balancing is as follows:

\begin{enumerate}
	\itemsep 0pt
	\item Initialize mesh (not distributed)
	\item (Re-)distribute the mesh (\routine{distr\_mesh} routine from \module{par\_mesh} module)
	\item Solve
	\item Estimate the error (end if solution is sufficiently accurate)
	\item Mark elements for refinement
	\item Refine the mesh ($hp$)
	\item Optionally, repartition the mesh for load balancing purposes (continue with Step 2); otherwise, continue with Step 3.
\end{enumerate}

In Step 6 (refinement), \hp3D implements subdomain inheritance when breaking nodes, i.e.~the subdomain value \var{subd} of a newly created node ($h$-refinement) is inherited from the father node. This may lead to \emph{load imbalance}. The \hp3D interface to the load balancing library Zoltan is provided in the module \module{zoltan\_wrapper}. This module includes several of the above mentioned functionalities to the user: for example, \routine{zoltan\_w\_eval} evaluates the quality of the current mesh, and \routine{zoltan\_w\_set\_lb(\var{LB})} sets an internal parameter (\var{ZOLTAN\_LB}) determining the load balancing strategy. The mesh repartitioning begins with a call to the routine \routine{zoltan\_w\_partition(\var{Mdle\_subd})} that returns, for each active element (i.e., active middle node) in the mesh, the newly assigned subdomain for the repartitioned mesh. Having determined a new partition for the mesh, the data migration step follows. This functionality is implemented in the \routine{distr\_mesh} routine of the module \module{par\_mesh}. 

Routines provided by the \module{zoltan\_wrapper} module:
\begin{itemize}
	\itemsep 0pt
	\item{\routine{zoltan\_w\_eval} \\
		Evaluates the quality of the current mesh partitioning.
	}
	\item{\routine{zoltan\_w\_set\_lb(\var{LB})} \\
		Sets global variable \var{ZOLTAN\_LB} to define partitioning algorithm.\\
		\var{ZOLTAN\_LB} $\in \{$BLOCK, RANDOM, RCB, RIB, HSFC, GRAPH$\}$.
	}
	\item{\routine{zoltan\_w\_partition(\var{Mdle\_subd(1:NRELES)})} \\
		Returns new partitioning (subdomain number for each middle node) for the active mesh.
	}
\end{itemize}

Each one of the partitioning algorithms is executed through the Zoltan library, however relies on partitioners from other external libraries (e.g.~ParMETIS); the list includes both geometry-based and graph-based algorithms. For more information, we refer to Zoltan's documentation \cite{ZoltanOverviewArticle2002}.

\paragraph{Data migration.} If during repartitioning element ownership has changed to a new subdomain, then the old owner and new owner of the element may exchange DOF data for each element node if requested by the user setting the global variable \code{\var{EXCHANGE\_DOF} = .true.} in module \module{par\_mesh}:
\begin{itemize}
	\itemsep -2pt
	\item For each element node, the old owner
	\begin{enumerate}
		\itemsep -2pt
		\item packs nodal data, i.e., $H^1$-, $H(\tcurl)$-, $H(\tdiv)$-, and $L^2$-DOFs, into a buffer array;
		\item sends nodal data to the new owner (point-to-point communication).
	\end{enumerate}
	\item For each element node, the new owner
	\begin{enumerate}
		\itemsep -2pt
		\item allocates nodal DOFs;
		\item receives nodal data from the old owner (point-to-point communication);
		\item unpacks nodal data from buffer array into $H^1$-, $H(\tcurl)$-, $H(\tdiv)$-, and $L^2$-DOFs.
	\end{enumerate}
\end{itemize}

Exchanging DOFs in this way is communication-intensive and should be deactivated by the user if not needed for the application. By default, \code{\var{EXCHANGE\_DOF} = .true.}.

%\input{Chapters/5_PARALLEL/comments}

%% file: Chapters/APPENDIX_A_EXAMPLES/appendix_examples_chapter.tex
%
%!TEX root = ../../hp3D_user_guide.tex
%

\chapter{Model Problems}
\label{chap:examples}

%--------------------------------------------------------------------

\input{Chapters/APPENDIX_A_EXAMPLES/sections/appendix_examples_poisson}

%%%%%%%%%%%%%%%%%%%%%%%%%%%%%%%

\input{Chapters/APPENDIX_A_EXAMPLES/sections/appendix_examples_elasticity}

%%%%%%%%%%%%%%%%%%%%%%%%%%%%%%%

\input{Chapters/APPENDIX_A_EXAMPLES/sections/appendix_examples_helmholtz}

%%%%%%%%%%%%%%%%%%%%%%%%%%%%%%%

\input{Chapters/APPENDIX_A_EXAMPLES/sections/appendix_examples_maxwell}

%%%%%%%%%%%%%%%%%%%%%%%%%%%%%%%

%\input{Chapters/APPENDIX_A_EXAMPLES/comments}

%% file: Chapters/APPENDIX_A_EXAMPLES/sections/appendix_examples_poisson.tex
%
%!TEX root = ../../../hp3D_user_guide.tex
%

%--------------------------------------------------------------------

\section{Poisson Problems}
\label{sec:poisson}

For the first model problem implementation, we consider the Poisson problem with inhomogeneous Dirichlet BC:
\begin{alignat*}{3}
	- \div \nabla u &= f && \quad \text{in } \Omega \, , \\
	u &= u_0 && \quad \text{on } \Gamma \, .
\end{alignat*}

Classical variational formulation:
\[
\left\{
\begin{array}{llll}
	u \in H^1(\Omega):  u = u_0 \text{ on } \Gamma \, , \\[5pt]
	(\nabla u, \nabla v) = (f,v) \, ,
	\quad v \in H^1(\Omega) \, : \, v = 0 \text{ on } \Gamma \, .
\end{array}
\right.
\]

%--------------------------------------------------------------------
\subsection{Galerkin implementation}
\label{sec:poisson-galerkin}

The Galerkin FE implementation of the variational problem is located in the application directory \file{problems/POISSON/GALERKIN}. In the remainder of this section, file paths may be given as relative paths within the application directory.

Input files:
\begin{itemize}
	\item{\file{control/control}: sets global control variables, e.g.
	\begin{itemize}
		\item \var{NEXACT} $\in \{$ \code{0,1} $\}$: indicates whether the exact solution is known.
		\item \var{EXGEOM} $\in \{$ \code{0,1} $\}$: indicates whether isoparametric or exact-geometry elements are used.
	\end{itemize}
	}
	\item{\file{input/physics}: sets initially allocated nodes and physics variables.
\begin{lstlisting}[caption=\file{POISSON/GALERKIN/input/physics} input file.]
100000              MAXNODS, nodes anticipated
1                   NR_PHYSA, physics attributes
field   contin  1   H1 variable
\end{lstlisting}
	\begin{itemize}
		\item The value of \var{MAXNODS} does not have to be precise; if more nodes are needed, the code allocates them on-the-fly. However, it is recommended for efficiency that the code does not reallocate, as well as not selecting \var{MAXNODS} much larger than needed.
		\item \code{\var{NR\_PHYSA}=1} specifies that \emph{one} physics variable is declared.
		\item \code{`field   contin  1'} specifies ``nickname, approximation space, number of components'' of a variable. The approximation spaces are: $H^1$ -- \code{contin}, $H(\tcurl)$ -- \code{tangen}, $H(\tdiv)$ -- \code{normal}, $L^2$ -- \code{discon}.
		\item For this Galerkin FE formulation, one $H^1$ variable is needed.
	\end{itemize}
	}
	\item{\file{geometries/hexa\_orient0}: defines the initial geometry mesh (a cube).}
\end{itemize}

Next, we take a look at the required routines that must be provided by the user:
\begin{itemize}
	\itemsep 0pt
	\item
	{\routine{set\_initial\_mesh}:\\
	for \emph{each} initial mesh element, this routine sets
	\begin{itemize}
		\itemsep 0pt
		\item the supported physics variables;
		\item the initial polynomial order of approximation;
		\item the boundary condition flags for element faces on the boundary.
	\end{itemize}
\begin{lstlisting}[caption=\file{POISSON/GALERKIN/}\routine{set\_initial\_mesh} routine.]
!..loop over initial mesh elements
   do iel=1,NRELIS

!  ...1. set physics
      ELEMS(iel)%nrphysics = 1
      allocate(ELEMS(iel)%physics(1))
      ELEMS(iel)%physics(1) ='field'

!  ...2. set initial order of approximation
      select case(ELEMS(iel)%etype)
         case(TETR); Nelem_order(iel) = 1*IP
         case(PYRA); Nelem_order(iel) = 1*IP
         case(PRIS); Nelem_order(iel) = 11*IP
         case(BRIC); Nelem_order(iel) = 111*IP
      end select

!  ...3. set BC flags: 0 - no BC ; 1 - Dirichlet; 2-9 Custom BCs
      ibc(1:6,1) = 0
      do ifc=1,nface(ELEMS(iel)%etype) ! loop through element faces
         neig = ELEMS(iel)%neig(ifc)
         select case(neig)
            case(0); ibc(ifc,1) = 1
         end select
      enddo

!  ...allocate BC flags (one per component), 
!     and encode face BCs into a single BC flag
      allocate(ELEMS(iel)%bcond(1))
      call encodg(ibc(1:6,1),10,6, ELEMS(iel)%bcond(1))
   enddo
\end{lstlisting}
	}
	\item
	{\routine{dirichlet}:
	\begin{itemize}
	\item User-provided routine required by the system routine \routine{update\_Ddof} which computes the Dirichlet DOFs for element nodes (vertices, edges, faces) with a non-homogeneous Dirichlet BCs. 
	\item \routine{update\_Ddof} interpolates $H^1$, $H(\tcurl)$, $H(\tdiv)$ Dirichlet data using \emph{projection-based} interpolation.\footnote{\fullcite{demkowicz2008interp}}
	\item Required only if non-homogeneous Dirichlet BCs were set by the user in \routine{set\_initial\_mesh}.
\end{itemize}

\begin{lstlisting}[caption=\file{POISSON/GALERKIN/common/}\routine{dirichlet} routine.]
!  routine dirichlet: returns Dirichlet data at a point
!   in:   Mdle          - middle node number
!         X             - a point in physical space
!         Icase         - node case (specifies supported variables)
!   out:  ValH, DvalH   - value of the H1 solution, 1st derivatives
!         ValE, DvalE   - value of the H(curl) solution, 1st derivatives
!         ValV, DvalV   - value of the H(div) solution, 1st derivatives
subroutine dirichlet(Mdle,X,Icase, ValH,DvalH,ValE,DvalE,ValV,DvalV)
\end{lstlisting}
	}
	\item
	{\routine{elem}:
	\begin{itemize}
	\item User-provided routine that computes the element-local stiffness matrix and load vector.
	\item System module \module{assembly} provides global arrays for this purpose:
	\begin{itemize}
		\item \var{ALOC(:,:)\%array}: Element-local stiffness matrix.
		\item \var{BLOC(:)\%array}: Element-local load vector.
		\item These arrays are declared \omp{omp threadprivate} for shared-memory parallel assembly of different element matrices with OpenMP threading.
	\end{itemize}
\end{itemize}

\routine{elem} is called during assembly for each middle node \var{Mdle} in the \emph{active mesh}.

\begin{remark}
Constrained approximation, modification for Dirichlet nodes, and static condensation of element-interior (bubble) DOFs are automatically done afterwards by the system routine \routine{celem\_system} which provides the \emph{modified element} matrices to the assembly procedure.
\end{remark}

\begin{lstlisting}[mathescape,caption=\file{POISSON/GALERKIN/}\routine{elem} routine]
!..determine element type; number of vertices, edges, and faces
   etype = NODES(Mdle)%ntype
   nrv = nvert(etype); nre = nedge(etype); nrf = nface(etype)
   
!..determine order of approximation
   call find_order(Mdle, norder)
   
!..determine edge and face orientations
   call find_orient(Mdle, norient_edge,norient_face)
   
!..determine nodes coordinates
   call nodcor(Mdle, xnod)
   
!..set quadrature points and weights
   call set_3D_int(etype,norder,norient_face, nrint,xiloc,waloc)

!  ....... element integrals:

!..loop over integration points
   do l=1,nrint

!  ...coordinates and weight of this integration point
      xi(1:3) = xiloc(1:3,l); wa = waloc(l)

!  ...H1 shape functions (for geometry)
      call shape3DH(etype,xi,norder,norient_edge,norient_face, nrdofH,shapH,gradH)

!  ...geometry map
      call geom3D(Mdle,xi,xnod,shapH,gradH,nrdofH, x,dxdxi,dxidx,rjac,iflag)

!  ...integration weight
      weight = rjac*wa

!  ...get the RHS
      call getf(Mdle,x, fval)

!  ...loop through H1 test functions
      do k1=1,nrdofH

!     ...Piola transformation: $q \rightarrow \hat q$ and $\nabla q \rightarrow J^{-T} \hat \nabla \hat q$
         q = shapH(k1)
         dq(1:3) = gradH(1,k1)*dxidx(1,1:3) + &
                    gradH(2,k1)*dxidx(2,1:3) + &
                    gradH(3,k1)*dxidx(3,1:3)

!     ...accumulate for the load vector: $(f,q)$
         b_loc(k1) = b_loc(k1) + q*fval*weight

!     ...loop through H1 trial functions
         do k2=1,nrdofH

!        ...Piola transformation: $p \rightarrow \hat p$ and $\nabla p \rightarrow J^{-T} \hat \nabla \hat p$
            p = shapH(k2)
            dp(1:3) = gradH(1,k2)*dxidx(1,1:3) + &
                       gradH(2,k2)*dxidx(2,1:3) + &
                       gradH(3,k2)*dxidx(3,1:3)

!        ...accumulate for the stiffness matrix: $(\nabla p, \nabla q)$
            a_loc(k1,k2) = a_loc(k1,k2) + weight*(dq(1)*dp(1)+dq(2)*dp(2)+dq(3)*dp(3))

   enddo; enddo; enddo
\end{lstlisting}
	}
\end{itemize}

This concludes the list of necessary input files and routines required for defining the application code from the library-perspective. However, the user is encouraged to take a look at the remaining files within the \file{POISSON/GALERKIN} directory which include the driver \routine{main} and a variety of auxiliary files. In a future version of the user manual, we may include a discussion of these auxiliary files as well.

%--------------------------------------------------------------------
\subsection{DPG primal implementation}
\label{sec:poisson-primal}

Broken primal DPG formulation:
\[
\left\{
\begin{array}{llll}
	(u, \hat \sigma_n) \in H^1(\Omega) \times H^{-1/2}(\Gamma_h): u = u_0 \text{ on } \Gamma \, , \\[5pt]
	(\nabla u, \nablah v) - \lb \hat \sigma_n , v \rb_{\Gamma_h} 
	= (f,v) \, ,\quad v \in H^1(\Omega_h) \, .
\end{array}
\right.
\]

Compared to the Galerkin FE implementation, the primal DPG implementation mostly differs in the \routine{elem} routine. In practice, the DPG method is implemented in its mixed form \cite{demkowicz2017dpg} but the extra unknown---the Riesz representation of the residual---is statically condensed on the element level (see Appendix~\ref{chap:dpg}). In the DPG implementation, the choice of the test norm is up to the user. This model problem employs the (broken) mathematicians norm as a test norm for the primal DPG implementation: $\| v \|^2_\test := (v,v) + (\nabla_h v, \nabla_h v)$.

The implementation is provided in \file{problems/POISSON/PRIMAL\_DPG}.

Compared to the input files for the Galerkin implementation, the only change is in the \file{physics} file:
\begin{lstlisting}[caption=\file{POISSON/PRIMAL\_DPG/input/physics} input file.]
100000              MAXNODS, nodes anticipated
2                   NR_PHYSA, physics attributes
field   contin  1   H1 variable
trace   normal  1   H(div) variable
\end{lstlisting}
We now have specified two physics unknowns---$u$ and $\hat \sigma_n$, i.e.~\code{\var{NR\_PHYSA}=2}. The additional trace unknown $\hat \sigma_n$ is declared as an $H(\tdiv)$ variable in the \file{physics} file; the normal trace $\hat \sigma_n$ must later be specified as such by setting \code{\var{PHYSAi(2)}=.true.} (e.g.~in the \routine{main} driver).

The \routine{elem} routine for DPG formulations can be structured into three distinct steps:
\vskip 5pt

\begin{minipage}{0.48\textwidth}
\begin{enumerate}
	\itemsep -10pt
	\item Element integration \vspace{-15pt}
	\begin{itemize}
		\itemsep -8pt
		\item Stiffness: $\mr B$
		\item Load: $\mr l$
		\item Gram matrix: $\mr G$
	\end{itemize}
	\item Boundary integration \vspace{-15pt}
	\begin{itemize}
		\itemsep -8pt
		\item Stiffness: $\mr{\hat B}$
	\end{itemize}
	\item Constructing DPG linear system \vspace{-15pt}
	\begin{itemize}
		\itemsep -8pt
		\item Dense linear algebra
		\item {Statically condensed system\\[-5pt] 
		stored in \var{ALOC}, \var{BLOC}}
	\end{itemize}
\end{enumerate}
\end{minipage}%
\begin{minipage}{0.48\textwidth}
\begin{figure}[H]
	\centering
	\begin{subfigure}[b]{0.6\textwidth}
		\includegraphics[height=1.8in]{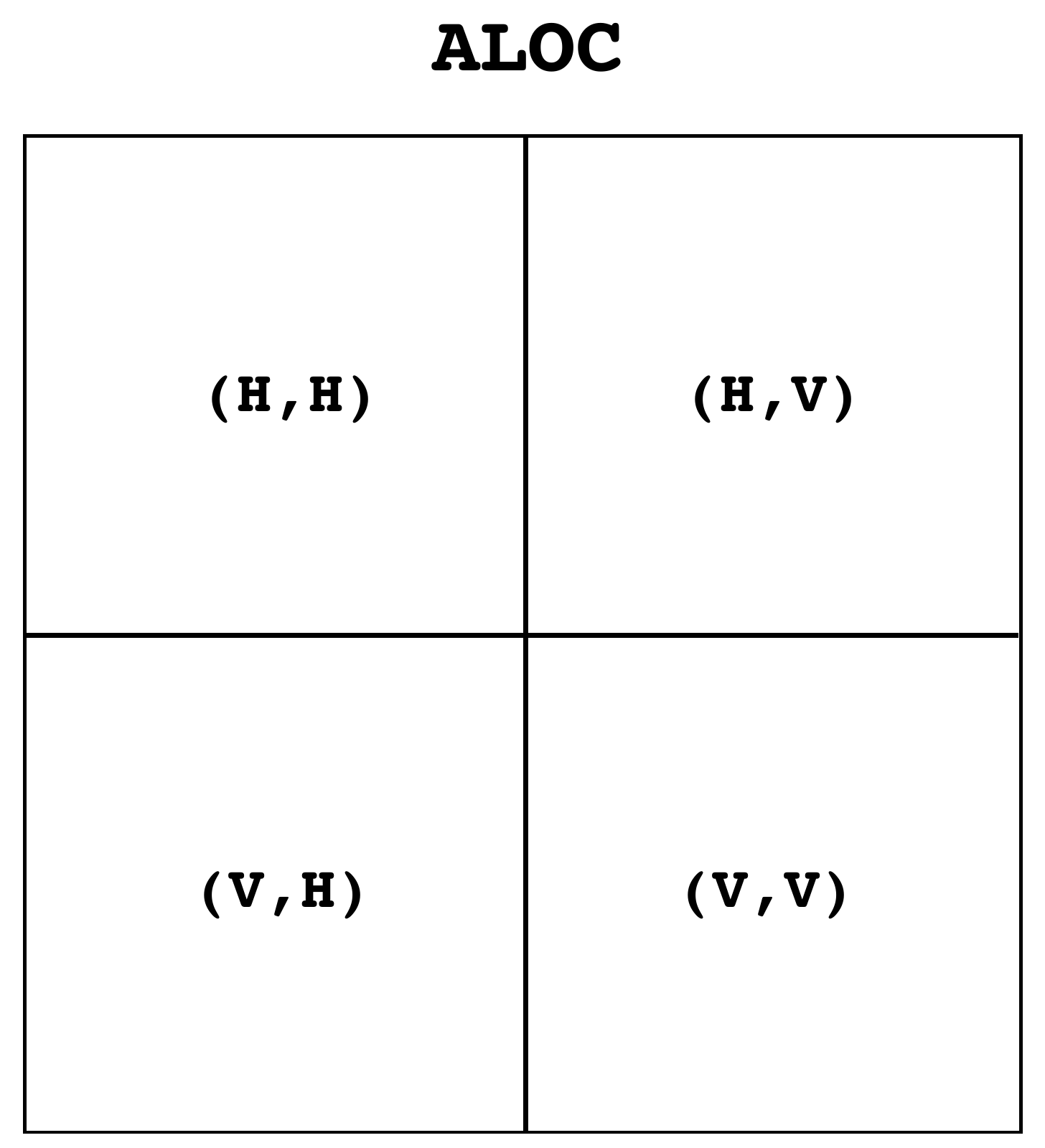}
	\end{subfigure}
	\begin{subfigure}[b]{0.3\textwidth}
		\includegraphics[height=1.8in]{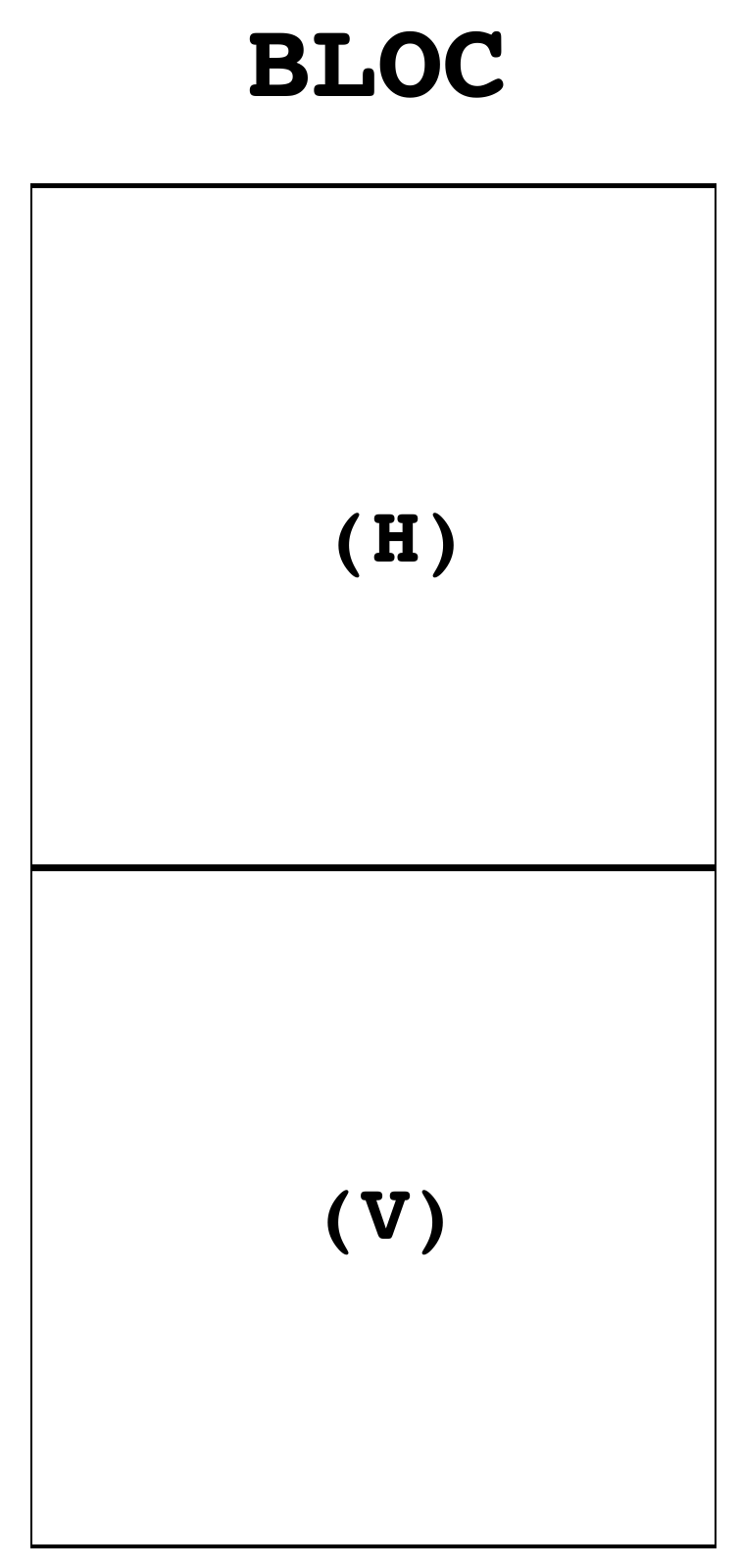}
	\end{subfigure}
	\caption*{Element-local system.}
\end{figure}
\end{minipage}

\vskip 5pt
\begin{minipage}[t]{0.60\textwidth}
Recall the statically condensed system\\[-5pt]
(cf.~Appendix~\ref{chap:dpg}): \vspace{-10pt}
\[
	\left[ \begin{array}{cc}
		\mr{B^* G^{-1} B} & \mr{B^* G^{-1} \hat B} \\
		\mr{\hat B^* G^{-1} B} & \mr{\hat B^* G^{-1} \hat B} \\
	\end{array} \right]
	\left[ \begin{array}{c}
		\mr{u_h} \\
		\mr{\hat u_h} 
	\end{array} \right]
	=
	\left[ \begin{array}{c}
		\mr{B^* G^{-1} l} \\
		\mr{\hat B^* G^{-1} l}
	\end{array} \right]
\]
\end{minipage}
\begin{minipage}[t]{0.39\textwidth}
Auxiliary local variables: \vspace{-10pt}
\begin{itemize}
	\itemsep -8pt
	\item \var{stiff\_HH} $\leftarrow \mr B$
	\item \var{stiff\_HV} $\leftarrow \mr{\hat B}$
	\item \var{bload\_H} $\leftarrow \mr l$
	\item \var{stiff\_ALL} $\leftarrow \left[ \mr{B \, | \, \hat B \, | \, l \, } \right]$
\end{itemize}
\end{minipage}

\vskip 5pt

In the \routine{elem} routine, these steps are implemented as follows:
\begin{enumerate}
	\item{ Preliminary set up:
\begin{lstlisting}[mathescape,caption=\file{POISSON/PRIMAL\_DPG/}\routine{elem}: preliminary set up]
!..allocate auxiliary matrices
   allocate(gramP(NrTest*(NrTest+1)/2))
   allocate(stiff_HH(NrTest,NrdofH))
   allocate(stiff_HV(NrTest,NrdofVi))
   allocate(bload_H(NrTest))

!..determine element type; number of vertices, edges, and faces
   etype = NODES(Mdle)%ntype
   nrv = nvert(etype); nre = nedge(etype); nrf = nface(etype)
   
!..determine order of approximation (element integrals)
   call find_order(Mdle, norder)
!..determine enriched order of approximation (hexa)
   nordP = NODES(Mdle)%order+NORD_ADD*111

!..determine edge and face orientations
   call find_orient(Mdle, norient_edge,norient_face)
!..determine nodes coordinates
   call nodcor(Mdle, xnod)
   
!  ....... element integrals
\end{lstlisting}
	}
	\item{ Element integration:
\begin{lstlisting}[mathescape,caption=\file{POISSON/PRIMAL\_DPG/}\routine{elem}: element integration]
!..use the enriched order to set the quadrature
   INTEGRATION = NORD_ADD ! $\Delta p \in \{1, 2, \ldots \}$
   call set_3D_int_DPG(etype,norder,norient_face, nrint,xiloc,waloc)

!..loop over integration points
   do l=1,nrint
!  ...coordinates and weight of this integration point
      xi(1:3) = xiloc(1:3,l); wa = waloc(l)

!  ...H1 shape functions (for geometry)
      call shape3DH(etype,xi,norder,norient_edge,norient_face, nrdofH,shapH,gradH)
!  ...discontinuous H1 shape functions
      call shape3HH(etype,xi,nordP, nrdof,shapHH,gradHH)

!  ...geometry map
      call geom3D(Mdle,xi,xnod,shapH,gradH,nrdofH, x,dxdxi,dxidx,rjac,iflag)
!  ...integration weight
      weight = rjac*wa
!  ...get the RHS
      call getf(Mdle,x, fval)

!  ...1st loop through enriched H1 test functions
      do k1=1,nrdofHH
!     ...Piola transformation
         v = shapHH(k1)
         dv(1:3) = gradHH(1,k1)*dxidx(1,1:3) + &
                    gradHH(2,k1)*dxidx(2,1:3) + &
                    gradHH(3,k1)*dxidx(3,1:3)
!
!     ...accumulate load: $(f,v)$
         bload_H(k1) = bload_H(k1) + fval*v*weight
!
!     ...loop through H1 trial functions
         do k2=1,nrdofH
!        ...Piola transformation
            dp(1:3) = gradH(1,k2)*dxidx(1,1:3) + &
                       gradH(2,k2)*dxidx(2,1:3) + &
                       gradH(3,k2)*dxidx(3,1:3)
!
!        ...accumulate stiffness: $(\nabla u, \nabla_h v)$
            stiff_HH(k1,k2) = stiff_HH(k1,k2) + weight*(dv(1)*dp(1)+dv(2)*dp(2)+dv(3)*dp(3))
         enddo

!     ...2nd loop through enriched H1 test functions for Gram matrix
         do k2=k1,nrdofHH
!        ...Piola transformation
            q = shapHH(k2)
            dq(1:3) = gradHH(1,k2)*dxidx(1,1:3) + & 
                       gradHH(2,k2)*dxidx(2,1:3) + &
                       gradHH(3,k2)*dxidx(3,1:3)

!        ...determine index in triangular packed format
            k = (k2-1)*k2/2+k1
!
!        ...accumulate Gram with test inner product: $(v,v)_\test := (v,v) + (\nabla_h v, \nabla_h v)$
            aux = q*v + (dq(1)*dv(1) + dq(2)*dv(2) + dq(3)*dv(3))
            gramP(k) = gramP(k) + aux*weight
         enddo; enddo; enddo
\end{lstlisting}
	}
	\item{ Boundary integration:
\begin{lstlisting}[mathescape,caption=\file{POISSON/PRIMAL\_DPG/}\routine{elem}: boundary integration.]
!..determine order of approximation (boundary integrals)
   norderi(1:nre+nrf) = norder(1:nre+nrf)
   norderi(nre+nrf+1) = 111

!..loop through element faces
   do ifc=1,nrf

!  ...sign factor to determine the outward normal unit vector
      nsign = nsign_param(etype,ifc)

!  ...face type ('tria','quad')
      ftype = face_type(etype,ifc)

!  ...face order of approximation
      call face_order(etype,ifc,norder, norderf)

!  ...set 2D quadrature
      INTEGRATION = NORD_ADD ! $\Delta p$
      call set_2D_int_DPG(ftype,norderf,norient_face(ifc), nrint,tloc,wtloc)

!  ...loop through integration points
      do l=1,nrint

!     ...face coordinates
         t(1:2) = tloc(1:2,l)
!     ...face parametrization
         call face_param(etype,ifc,t, xi,dxidt)

!     ...determine discontinuous H1 shape functions
         call shape3HH(etype,xi,nordP, nrdof,shapHH,gradHH)
!     ...determine element H(div) shape functions (for fluxes), interfaces only (no bubbles)
         call shape3DV(etype,xi,norderi,norient_face, nrdof,shapV,divV)

!     ...determine element H1 shape functions (for geometry)
         call shape3DH(etype,xi,norder,norient_edge,norient_face, nrdof,shapH,gradH)
!     ...geometry map
         call bgeom3D(Mdle,xi,xnod,shapH,gradH,nrdofH,dxidt,nsign, &
                        x,dxdxi,dxidx,rjac,dxdt,rn,bjac)
!     ...integration weight
         weight = bjac*wtloc(l)

!     ...loop through enriched H1 test functions
         do k1=1,nrdofHH
            v = shapHH(k1)

!        ...loop through H(div) trial functions
            do k2=1,nrdofVi
!           ...Piola transformation
               s(1:3) = (dxdxi(1:3,1)*shapV(1,k2) + &
                          dxdxi(1:3,2)*shapV(2,k2) + &
                          dxdxi(1:3,3)*shapV(3,k2)) / rjac
!           ...normal component
               sn = s(1)*rn(1)+s(2)*rn(2)+s(3)*rn(3)
!
!           ...accumulate stiffness: $-\lb \sigma \cdot n, v \rb_{\Gamma_h}$
               stiff_HV(k1,k2) = stiff_HV(k1,k2) - sn*v*weight

            enddo; enddo ! end loop through trial / test functions
   enddo; enddo; ! end loop through integration points / faces
\end{lstlisting}
	}
	\item{ Construction of DPG linear system:
\begin{lstlisting}[mathescape,caption=\file{POISSON/PRIMAL\_DPG/}\routine{elem}: constructing DPG linear system.]
!---------------------------------------------------------------------
!  Construction of statically condensed DPG linear system
!---------------------------------------------------------------------

!..create auxiliary matrix for dense linear algebra
   allocate(stiff_ALL(NrTest,NrTrial+1))

!..Total test/trial DOFs of the element
   i = NrTest ; j1 = NrdofH ; j2 = NrdofVi

!..Copy stiffness and load into one matrix: $\text{\var{stiff\_ALL}} \leftarrow [ \mr{B \, | \, \hat B \, | \, l \, } ]$
   stiff_ALL(1:i,1:j1)          = stiff_HH(1:i,1:j1)
   stiff_ALL(1:i,j1+1:j1+j2) = stiff_HV(1:i,1:j2)
   stiff_ALL(1:i,j1+j2+1)       = bload_H(1:i)

   deallocate(stiff_HH,stiff_HV)

!..A. Compute Cholesky factorization of Gram Matrix, $\mr{G=U^T U \, (=LL^T)}$
   call DPPTRF('U',NrTest,gramP,info)
!
!..B. Solve triangular system to obtain $\mr{\tilde B,\ \text{i.e.~solve } (L \tilde B=)\, U^T \tilde B = [B \, | \, l]}$
   call DTPTRS('U','T','N',NrTest,NrTrial+1,gramP,stiff_ALL,NrTest,info)

   allocate(raloc(NrTrial+1,NrTrial+1)); raloc = ZERO

!..C. Matrix multiply: $\mr{B^T G^{-1} B \, (=\tilde B^T \tilde B)}$
   call DSYRK('U','T',NrTrial+1,NrTest,ZONE,stiff_ALL,NrTest,ZERO,raloc,NrTrial+1)

!..D. Fill lower triangular part of Hermitian matrix $\mr{\tilde B^T \tilde B}$
   do i=1,NrTrial
      raloc(i+1:NrTrial+1,i) = raloc(i,i+1:NrTrial+1)
   enddo

!..raloc has now all blocks of the stiffness and load:
!  $r_{\text{aloc}} = \var{ALOC(1,1)} \quad \var{ALOC(1,2)} \quad \var{BLOC(1)}$
!         $\, \var{ALOC(2,1)} \quad \var{ALOC(2,2)} \quad \var{BLOC(2)}$
\end{lstlisting}
	}
\end{enumerate}

%--------------------------------------------------------------------
\input{Chapters/APPENDIX_A_EXAMPLES/sections/appendix_examples_UW_poisson}

%% file: Chapters/APPENDIX_A_EXAMPLES/sections/appendix_examples_UW_poisson.tex
%
%!TEX root = ../../../hp3D_user_guide.tex
%

%--------------------------------------------------------------------

\subsection{DPG ultraweak implementation}
\label{sec:poisson-ultraweak}

The ultraweak formulation is based on the first-order system:
\begin{align*}
	\bs \sigma - \grad u &= \bs 0 , \\
	- \div \bs \sigma &= f .
\end{align*}

The broken ultraweak formulation is then given by:
\[
\left\{
\begin{array}{llll}
	\hat u \in H^{1/2}(\Gamma_h): \hat u = u_0 \text{ on } \Gamma,\,
	\hat{\sigma}_n \in H^{-1/2}(\Gamma_h),\,
	u \in L^2(\Omega),\,
	\bs \sigma \in \bs{L^2}(\Omega) , \\[5pt]
	(\bs \sigma, \bs \tau) + (u, \nablah \cdot \bs \tau)
	- \lb \hat u, \bs \tau \cdot \bs n \rb_{\Gamma_h}
	+ (\bs \sigma, \nablah v)
	- \lb \hat{\sigma}_n , v \rb_{\Gamma_h}
	= (f,v) \, , \
	v \in H^1(\Omega_h), 
	\bs \tau \in H(\tdiv; \Omega_h) \, .
\end{array}
\right.
\]

The DPG ultraweak implementation differs from the primal DPG implementation mostly in the \routine{elem} routine. As before, the static condensation of the Riesz representation of the residual is done at the element level (cf.~Section~\ref{sec:poisson-primal}). This model problem employs the adjoint graph norm as a test norm for the ultraweak DPG implementation:
\[
	\| (v,\bs \tau) \|^2_\test := 
	\| \nablah v + \bs \tau \|^2 + \| \nablah \cdot \bs \tau \|^2
	+ \| v \|^2 + \| \bs \tau \|^2 .
\]

The implementation is provided in \file{problems/POISSON/ULTRAWEAK\_DPG}.

\begin{lstlisting}[caption=\file{POISSON/ULTRAWEAK\_DPG/input/physics} input file.]
100000               MAXNODS, nodes anticipated
4                    NR_PHYSA, physics attributes
trace_a contin  1    H1 variable
trace_b normal  1    H(div) variable
field    discon  1   L2 variable
grad     discon  3   L2 Variable
\end{lstlisting}

In the ultraweak formulation, there are four physics unknowns: $\hat u$, $\hat{\sigma}_n$, $u$, $\bs \sigma$; i.e.~\code{\var{NR\_PHYSA}=4}. The variables are defined in the \file{physics} file and the continuous trace $\hat u$ and the normal trace $\hat \sigma_n$ are specified as such by setting \code{\var{PHYSAi(1:2)}=.true.}. As before, the DPG \routine{elem} routine can be split into three distinct steps:

\begin{minipage}{0.48\textwidth}
\begin{enumerate}
	\itemsep -10pt
	\item Element integration \vspace{-15pt}
	\begin{itemize}
		\itemsep -8pt
		\item Stiffness: $\mr B$
		\item Load: $\mr l$
		\item Gram matrix: $\mr G$
	\end{itemize}
	\item Boundary integration \vspace{-15pt}
	\begin{itemize}
		\itemsep -8pt
		\item Stiffness: $\mr{\hat B_{1,1}, \hat B_{1,2}, \hat B_{2,1}, \hat B_{2,2}}$
	\end{itemize}
	\item Constructing DPG linear system \vspace{-15pt}
	\begin{itemize}
		\itemsep -8pt
		\item Dense linear algebra
		\item {Statically condensed system\\[-5pt] 
		stored in \var{ALOC}, \var{BLOC}}
	\end{itemize}
\end{enumerate}
\end{minipage}%
\begin{minipage}{0.48\textwidth}
\begin{figure}[H]
	\centering
		\includegraphics[scale=0.5]{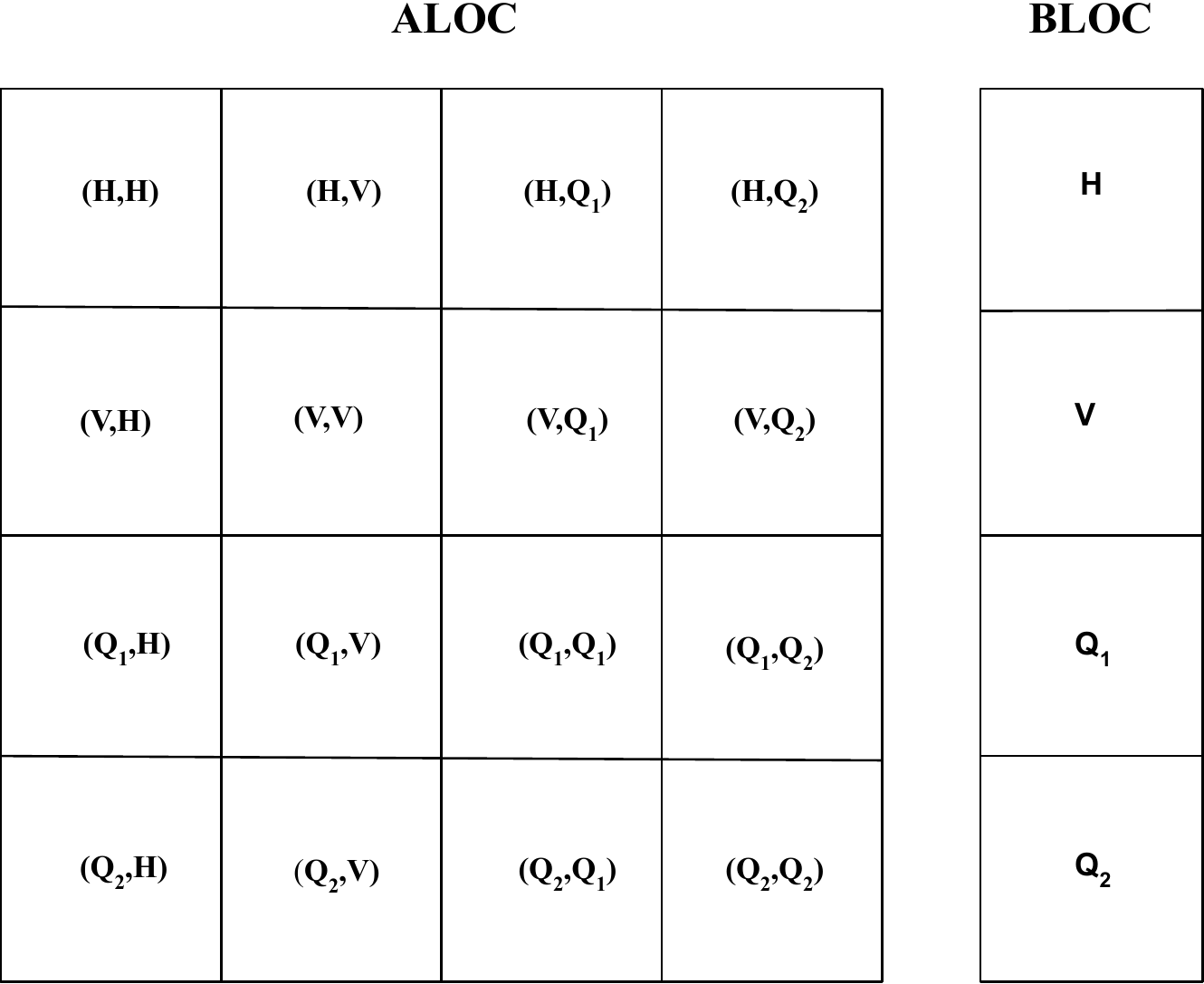}
	\caption*{Element-local system.}
\end{figure}
\end{minipage}

\vskip 5pt
%\begin{minipage}[t]{0.60\textwidth}
%%Recall the statically condensed system\\[-5pt]
%(cf.~Appendix~\ref{chap:dpg}): \vspace{-10pt}
%\[
%  \var{stiff\_ALL}^* G^{-1} \var{stiff\_ALL} = \var{stiff\_ALL}^* G^{-1} \var{bload\_H} 
%\]
%\end{minipage}
\begin{minipage}[t]{\textwidth}
Auxiliary local variables: \vspace{-10pt}
\begin{itemize}
	\itemsep -8pt
	\item \var{stiff\_HH} $\leftarrow \mr{\hat B_{1,1}}$
	\item \var{stiff\_HV} $\leftarrow \mr{\hat B_{1,2}}$
	\item \var{stiff\_HQ1} $\leftarrow \mr{B_{1,3}}$
	\item \var{stiff\_HQ2} $\leftarrow \mr{B_{1,4}}$
	\item \var{stiff\_VH} $\leftarrow \mr {\hat B_{2,1}}$
	\item \var{stiff\_VV} $\leftarrow \mr{\hat B_{2,2}}$
	\item \var{stiff\_VQ1} $\leftarrow \mr{B_{2,3}}$
	\item \var{stiff\_VQ2} $\leftarrow \mr{B_{2,4}}$
	
	\item \var{bload\_H} $\leftarrow \mr {\left[ {l} \, | \, 0\right]}^T$
	\item \vskip 4pt \var{stiff\_ALL} $\leftarrow \left[\begin{array}{c|c|c|c|c}
	\hat B_{1,1} \, | \, & \hat B_{1,2} \, | \, & B_{1,3} \, | \, & B_{1,4} \, | \, & \text{l} \\ \hline
	\hat B_{2,1} \, | \, & \hat B_{2,2} \, | \, & B_{2,3} \, | \, & B_{2,4} \, | \, & 0
	\end{array} \right]$
\end{itemize}
\end{minipage}

The preliminary setup for the ultraweak formulation is similar to the setup of the primal DPG implementation with additional memory allocation for additional unknowns and test functions. The element integration routine for the ultraweak DPG formulation is given by the following code:
\begin{lstlisting}[mathescape,caption=\file{POISSON/ULTRAWEAK\_DPG/}\routine{elem}: element integration]
!..use the enriched order to set the quadrature
   INTEGRATION = NORD_ADD ! $\triangle p \in \{1,2,...\}$
   call set_3D_int_DPG(etype,norder,norient_face, nrint,xiloc,waloc)

!..loop over integration points
   do l=1,nrint

!..coordinates and weight of this integration point
      xi(1:3) = xiloc(1:3,l)
      wa = waloc(l)

!..H1 shape functions (for geometry)
      call shape3DH(etype,xi,norder,norient_edge,norient_face, nrdof,shapH,gradH)

!..L2 shape functions
      call shape3DQ(etype,xi,norder, nrdof,shapQ)

!..discontinuous H1 shape functions
      call shape3HH(etype,xi,nordP, nrdof,shapHH,gradHH)
!..discontinuous H(div) shape functions
      call shape3VV(etype,xi,nordP, nrdof,shapVV,divVV)

!..geometry map
      call geom3D(Mdle,xi,xnod,shapH,gradH,NrdofH, x,dxdxi,dxidx,rjac,iflag)

!..integration weight
      weight = rjac*wa

!..get the RHS
      call getf(Mdle,x, fval)

!..1st loop through enriched H1 test functions
      do k1=1,NrdofHH
!..Piola transformation
         v = shapHH(k1)
         dv(1:3) = gradHH(1,k1)*dxidx(1,1:3) &
                  + gradHH(2,k1)*dxidx(2,1:3) &
                  + gradHH(3,k1)*dxidx(3,1:3)

!..accumulate load vector: $\mr{(f,v)}$
         bload_H(k1) = bload_H(k1) + fval*v*weight

!..loop through L2 trial functions
         do k2=1,NrdofQ
!..loop over components of $\sigma$
            do ivar = 1,3
               n1 = k1; n2 = (k2 - 1)*3 + ivar
               sig = ZERO
!..Piola Transform for $\mr{{ivar}^{th}}$ L2 component of $\sigma$
               sig(ivar) = shapQ(k2) / rjac
!..accumulate stiffness: $\mr{(\bs \sigma,\nabla v)}$
               stiff_HQ2(n1,n2) = stiff_HQ2(n1,n2) + weight * &
               	                 (sig(1)*dv(1) + sig(2)*dv(2) + sig(3)*dv(3))
            enddo; enddo !..end of loop through L2 trial functions

!..2nd loop through enriched H1 test functions for Gram matrix
         do k2=k1,NrdofHH
!..Piola transformation
            q = shapHH(k2)
            dq(1:3) = gradHH(1,k2)*dxidx(1,1:3) &
                     + gradHH(2,k2)*dxidx(2,1:3) &
                     + gradHH(3,k2)*dxidx(3,1:3)

!..determine index in triangular packed format            
            k = nk(k1,k2)
!..accumulate components of Gram matrix corresponding to $\mr{v}$
            aux = q*v + (dq(1)*dv(1) + dq(2)*dv(2) + dq(3)*dv(3))
            gramP(k) = gramP(k) + aux*weight
         enddo !..end of 2nd loop through enriched H1 test functions

!..cross terms for  graph norm
         do k2 = 1,NrdofVV
            tau_a(1:3) = ( dxdxi(1:3,1) * shapVV(1,k2)       &
                           + dxdxi(1:3,2) * shapVV(2,k2)     &
                           + dxdxi(1:3,3) * shapVV(3,k2) ) / rjac

            k = nk(k1,NrdofHH+k2)
!..accumulate components of Gram matrix corresponding to inner products involving $\mr{v}$ and $\mr{\bs \tau}$
            aux = dv(1) * tau_a(1) + dv(2) * tau_a(2) + dv(3) * tau_a(3)
            gramP(k) = gramP(k) + aux * weight
         enddo
      enddo !..end of 1st loop through enriched H1 test functions

!..loop over discontinuous H(div) test functions
      do k1=1,NrdofVV
!..Piola transformation
         divtau_a = divVV(k1)/rjac
         tau_a(1:3) = ( dxdxi(1:3,1) * shapVV(1,k1)      &
                       + dxdxi(1:3,2) * shapVV(2,k1)     &
                       + dxdxi(1:3,3) * shapVV(3,k1) ) / rjac
         
         do k2=1,NrdofQ
            u = shapQ(k2) / rjac !..Piola transformation

!..accumulate stiffness: $\mr{(u,\nabla \cdot \bs \tau)}$
            stiff_VQ1(k1,k2) = stiff_VQ1(k1,k2) + weight*(u * divtau_a)

            do ivar = 1,3
               sig = ZERO
               sig(ivar) = shapQ(k2)/rjac  ! Piola transform for the $\mr{{ivar}^{th}}$ component
               n1 = k1; n2 = (k2 - 1)*3 + ivar
!..accumulate stiffness: $\mr{(\bs \sigma,\bs \tau)}$
               stiff_VQ2(n1,n2) = stiff_VQ2(n1,n2) + weight * & 
                                  (tau_a(1)*sig(1) + tau_a(2)*sig(2) + tau_a(3)*sig(3))
            enddo; enddo ! end of loop over L2 trial functions

!..Gram matrix contribution for H(div) inner product
         do k2=k1,NrdofVV
!..Piola transformation
            divtau_b = divVV(k2)/rjac
            tau_b(1:3) = ( dxdxi(1:3,1) * shapVV(1,k2)       &
                           + dxdxi(1:3,2) * shapVV(2,k2)     &
                           + dxdxi(1:3,3) * shapVV(3,k2) ) / rjac
!
            k = nk(k1 + NrdofHH,k2 + NrdofHH)
            aux = divtau_a * divtau_b + 2.d0 * 
                  (tau_a(1)*tau_b(1) + tau_a(2)*tau_b(2) + tau_a(3)*tau_b(3))
            gramP(k) = gramP(k) +  weight * aux
         enddo; enddo !..end of loop over H(div) test functions
   enddo !..end of loop over integration points
\end{lstlisting}

Next, we provide the code performing the boundary integration:
\begin{lstlisting}[mathescape,caption=\file{POISSON/ULTRAWEAK\_DPG/}\routine{elem}: boundary integration]
!..loop through element faces
   do ifc=1,nrf

!..sign factor to determine the outward normal unit vector
      nsign = nsign_param(etype,ifc)

!..face type
      ftype = face_type(etype,ifc)

!..face order of approximation
      call face_order(etype,ifc,norder, norderf)

!..set 2D quadrature
      INTEGRATION = NORD_ADD
      call set_2D_int_DPG(ftype,norderf,norient_face(ifc), nrint,tloc,wtloc)

!..loop through integration points
      do l=1,nrint

!..face coordinates
         t(1:2) = tloc(1:2,l)

!..face parametrization
         call face_param(etype,ifc,t, xi,dxidt)

!..determine discontinuous H1 shape functions
         call shape3HH(etype,xi,nordP, nrdof,shapHH,gradHH)
!..discontinuous H(div) shape functions
         call shape3VV(etype,xi,nordP, nrdof,shapVV,divVV)

!..determine element H1 shape functions (for geometry)
         call shape3DH(etype,xi,norder,norient_edge,norient_face, nrdof,shapH,gradH)

!..determine H(div) shape functions (for fluxes) on face
         call shape3DV(etype,xi,norderi,norient_face, nrdof,shapV,divV)

!..geometry map
         call bgeom3D(Mdle,xi,xnod,shapH,gradH,NrdofH,dxidt,nsign, &
                       x,dxdxi,dxidx,rjac,dxdt,rn,bjac)
!..integration weight
         weight = bjac*wtloc(l)

!..loop through enriched H1 test functions
         do k1=1,NrdofHH
            v = shapHH(k1)

!..loop through H(div) trial functions on face
            do k2=1,NrdofVi_b
!..Piola transformation
               s(1:3) = ( dxdxi(1:3,1)*shapV(1,k2) &
                         + dxdxi(1:3,2)*shapV(2,k2) &
                         + dxdxi(1:3,3)*shapV(3,k2) ) / rjac
!..normal component
               sn = s(1)*rn(1) + s(2)*rn(2) + s(3)*rn(3)

!..accumulate stiffness: $\mr{- \left\langle \hat{\sigma}_n,v \right\rangle}$
               stiff_HV(k1,k2) = stiff_HV(k1,k2) - sn*v*weight
            enddo !..end of loop through H(div) trial functions on face
         enddo !..end of loop through enriched H1 test functions

!..loop through enriched H(div) test functions
         do k1=1,NrdofVV
!..Piola Transform
            tau_a(1:3) = ( dxdxi(1:3,1) * shapVV(1,k1)      &
                          + dxdxi(1:3,2) * shapVV(2,k1)     &
                          + dxdxi(1:3,3) * shapVV(3,k1) ) / rjac
            tn = tau_a(1)*rn(1) + tau_a(2)*rn(2) + tau_a(3)*rn(3)

            do k2=1,NrdofVi_a
!..accumulate stiffness: $\mr{-\left\langle \hat u, \tau \cdot n \right\rangle}$
               u_hat = shapH(k2)
               stiff_VH(k1,k2) = stiff_VH(k1,k2) - weight * tn * u_hat
            enddo !..end of loop through H1 trial functions
         enddo !..end of loop through enriched H(div) test functions
      enddo !..end of loop through integration points
   enddo !..end of loop through element faces
\end{lstlisting}

Finally, the code performing the construction of the DPG linear system:
\begin{lstlisting}[mathescape,caption=\file{POISSON/ULTRAWEAK\_DPG/}\routine{elem}: constructing DPG linear system.]
   allocate(stiff_ALL(NrTest,NrTrial+1))
   stiff_ALL = ZERO

!  Total test/trial DOFs of the element
   i1 = NrdofHH; i2 = NrdofVV; j1 = NrdofVi_a; j2 = NrdofVi_b; j3 = NrdofU; j4 = NrdofS

!..Copy stiffness and load into one matrix
   stiff_ALL(1:i1,1:j1) = stiff_HH
   stiff_ALL(1:i1,j1+1:j1+j2) = stiff_HV
   stiff_ALL(1:i1,j1+j2+1:j1+j2+j3) = stiff_HQ1
   stiff_ALL(1:i1,j1+j2+j3+1:j1+j2+j3+j4) = stiff_HQ2

   stiff_ALL(i1+1:i1+i2,1:j1) = stiff_VH
   stiff_ALL(i1+1:i1+i2,j1+1:j1+j2) = stiff_VV
   stiff_ALL(i1+1:i1+i2,j1+j2+1:j1+j2+j3) = stiff_VQ1
   stiff_ALL(i1+1:i1+i2,j1+j2+j3+1:j1+j2+j3+j4) = stiff_VQ2

   stiff_ALL(1:i1+i2,j1+j2+j3+j4+1) = bload_H

!..A. Compute Cholesky factorization of Gram Matrix, $\mr{G=U^T U (=LL^T)}$
   call DPPTRF('U',NrTest,gramP,info)

!..B. Solve triangular system to obtain $\mr{\tilde{B}}$, $\mr{(LX=) U^T X = [B|l]}$
   call DTPTRS('U','T','N',NrTest,NrTrial+1,gramP,stiff_ALL,NrTest,info)

   allocate(raloc(NrTrial+1,NrTrial+1)); raloc = ZERO

!..C. Matrix multiply: $\mr{B^T G^{-1} B (=\tilde{B}^T \tilde{B})}$
   call DSYRK('U','T',NrTrial+1,NrTest,ZONE,stiff_ALL,NrTest,ZERO,raloc,NrTrial+1)

!..D. Fill lower triangular part of Hermitian matrix using the upper triangular matrix $\mr{\tilde B^T \tilde B}$
   do i=1,NrTrial
      raloc(i+1:NrTrial+1,i) = raloc(i,i+1:NrTrial+1)
   enddo
   
!..raloc has now all blocks of the stiffness and load:
!  $r_{\text{aloc}} =   \var{ALOC(1,1)} \quad \var{ALOC(1,2)} \quad \var{ALOC(1,3)} \quad \var{ALOC(1,4)} \quad\var{BLOC(1)}$
!         $\, \var{ALOC(2,1)} \quad \var{ALOC(2,2)} \quad \var{ALOC(2,3)} \quad \var{ALOC(2,4)} \quad \var{BLOC(2)}$
!         $\, \var{ALOC(3,1)} \quad \var{ALOC(3,2)} \quad \var{ALOC(3,3)} \quad \var{ALOC(3,4)} \quad \var{BLOC(3)}$
!         $\, \var{ALOC(4,1)} \quad \var{ALOC(4,2)} \quad \var{ALOC(4,3)} \quad \var{ALOC(4,4)} \quad \var{BLOC(4)}$
\end{lstlisting}

%% file: Chapters/APPENDIX_A_EXAMPLES/sections/appendix_examples_elasticity.tex
%
%!TEX root = ../../../hp3D_user_guide.tex
%

%--------------------------------------------------------------------

\section{Linear Elasticity Problems}
\label{sec:elasticity}

%\subsection{Galerkin implementation}
%\label{subsec:elasticity-galerkin}

Linear elasticity will be added in a future version of the user manual.

%The general structure of the \file{ELASTICITY} application directory is similar to the implementation of previously discussed applications. This section only focuses on the specific aspects of implementing particular formulations of linear elasticity. We encourage reading Section~\ref{sec:poisson} on Poisson problems for additional discussion of problem implementations, in general.

%% file: Chapters/APPENDIX_A_EXAMPLES/sections/appendix_examples_helmholtz.tex
%
%!TEX root = ../../../hp3D_user_guide.tex
%

%--------------------------------------------------------------------

\section{Helmholtz Problems}
\label{sec:helmholtz}

%\subsection{Galerkin implementation}
%\label{subsec:helmholtz-galerkin}

Helmholtz (linear acoustics) will be added in a future version of the user manual.

%The general structure of the \file{HELMHOLTZ} application directory is similar to the implementation of previously discussed applications. This section only focuses on the specific aspects of implementing particular formulations of the Helmholtz problem. We encourage reading Section~\ref{sec:poisson} on Poisson problems for additional discussion of problem implementations, in general.

%% file: Chapters/APPENDIX_A_EXAMPLES/sections/appendix_examples_maxwell.tex
%
%!TEX root = ../../../hp3D_user_guide.tex
%

%--------------------------------------------------------------------

\section{Maxwell Problems}
\label{sec:maxwell}

For time-harmonic Maxwell problems, the solution is complex-valued; the \hp3D library must therefore be compiled with preprocessing flag \code{\var{COMPLEX}=1}.

The general structure of the \file{MAXWELL} application directory is similar to the implementation of previously discussed applications. This section only focuses on the specific aspects of implementing particular formulations of the Maxwell problem. We encourage reading Section~\ref{sec:poisson} on Poisson problems for additional discussion of problem implementations, in general.

\begin{itemize}
\item
{
Linear time-harmonic Maxwell equations:
\begin{alignat*}{3}
	\curl \bs E + i \omega \mu \bs H
	&= \bs 0 &\quad \text{in } \Omega \, , \\
	\curl \bs H - (i \omega \eps + \sigma) \bs E 
	&= \bs J^{\text{imp}} &\quad \text{in } \Omega \, , \\
	\bs n \times \bs E &= \bs n \times \bs E_0 &\quad \text{on } \Gamma \, .
\end{alignat*}
}
\item
{
Curl--curl formulation:
\begin{alignat*}{3}
	\curl (\mu^{-1} \curl \bs E) - (\omega^2 \eps - i \omega \sigma) \bs E
	&= -i \omega \bs J^{\text{imp}}  &\quad \text{in } \Omega \, , \\
	\bs n \times \bs E &= \bs n \times \bs E_0 &\quad \text{on } \Gamma \, .
\end{alignat*}
}
\item
{
Classical variational formulation:
\[
\left\{
\begin{array}{lll}
	\bs E \in \Hcurl : \bs n \times \bs E = \bs n \times \bs E_0 \text{ on } \Gamma \, , \\[5pt]
	(\mu^{-1} \curl \bs E,\curl \bs F) - ((\omega^2 \eps - i \omega \sigma) \bs E, \bs F)
	= -i \omega (\bs J^{\text{imp}}, \bs F) \, , \\[5pt]
	\hfill
	\quad \bs F \in \Hcurl \, : \, \bs n \times \bs F = \bs 0 \text{ on } \Gamma \, .
\end{array}
\right.
\]
The formulation involves just one unknown $\bs E \in \Hcurl$, defined on the whole domain.
}
\end{itemize}

\subsection{Galerkin implementation}
\label{subsec:maxwell-galerkin}

The Galerkin FE implementation of the classical (curl-curl) variational problem is located in the application directory \file{problems/MAXWELL/GALERKIN}. In the remainder of this section, file paths may be given as relative paths within the application directory.

The physics input file defines the vector-valued unknown: electric field $\bs E \in \Hcurl$.
\begin{lstlisting}[caption=\file{MAXWELL/GALERKIN/input/physics} input file.]
100000              MAXNODS, nodes anticipated
1                   NR_PHYSA, physics attributes
field  tangen  1    H(curl) variable
\end{lstlisting}

The curl-curl formulation above specifies a non-homogeneous Dirichlet boundary condition for the tangential component of the electric field. When setting a Dirichlet flag for an $H(\tcurl)$-variable in the \routine{set\_initial\_mesh} routine, the code automatically interprets the BC to be applied to the tangential component. As in previous examples, the values of the field and its derivative on the boundary must be supplied by the user in the \routine{dirichlet} routine.

The \routine{elem} routine for the Galerkin implementation of the Maxwell problem is given by the following code:

\begin{lstlisting}[mathescape,caption=\file{MAXWELL/GALERKIN/}\routine{elem} routine]
!..determine element type
   etype = NODES(Mdle)%ntype
   
!..determine order of approximation
   call find_order(Mdle, norder)
   
!..determine edge and face orientations
   call find_orient(Mdle, norient_edge,norient_face)
   
!..determine nodes coordinates
   call nodcor(Mdle, xnod)
   
!..set quadrature points and weights
   call set_3D_int(etype,norder,norient_face, nrint,xiloc,waloc)

!  ....... element integrals:

!..loop over integration points
   do l=1,nrint

!  ...coordinates and weight of this integration point
      xi(1:3) = xiloc(1:3,l); wa = waloc(l)

!  ...H1 shape functions (for geometry)
      call shape3DH(etype,xi,norder,norient_edge,norient_face, nrdofH,shapH,gradH)

!  ...H(curl) shape functions
      call shape3DE(etype,xi,norder,norient_edge,norient_face, nrdofE,shapE,curlE)

!  ...geometry map
      call geom3D(Mdle,xi,xnod,shapH,gradH,nrdofH, x,dxdxi,dxidx,rjac,iflag)

!  ...integration weight
      weight = rjac*wa

!  ...get the RHS (complex, vector-valued)
      call getf(Mdle,x, zJ)

!  ...loop through H(curl) test functions
      do k1=1,nrdofE

!     ...Piola transformation: $F \rightarrow J^{-T} \hat F$ and $\curl F \rightarrow (\det J)^{-1} J \, \hat \nabla \times \hat F$
         F(1:3) = shapE(1,k1)*dxidx(1,1:3) &
                 + shapE(2,k1)*dxidx(2,1:3) &
                 + shapE(3,k1)*dxidx(3,1:3)
         CF(1:3) = ( dxdxi(1:3,1)*curlE(1,k1) &
                    + dxdxi(1:3,2)*curlE(2,k1) &
                    + dxdxi(1:3,3)*curlE(3,k1) ) / rjac

!     ...accumulate for the load vector: $(-i \omega J,F)$
         za = F(1)*zJ(1) + F(2)*zJ(2) + F(3)*zJ(3)
         Zbloc(k1) = Zbloc(k1) - ZI*OMEGA*za*weight

!     ...loop through H(curl) trial functions
         do k2=1,nrdofE

!        ...Piola transformation: $E \rightarrow J^{-T} \hat E$ and $\curl E \rightarrow (\det J)^{-1} J \, \hat \nabla \times \hat E$
            E(1:3) = shapE(1,k2)*dxidx(1,1:3) &
                    + shapE(2,k2)*dxidx(2,1:3) &
                    + shapE(3,k2)*dxidx(3,1:3)
            CE(1:3) = ( dxdxi(1:3,1)*curlE(1,k2) &
                       + dxdxi(1:3,2)*curlE(2,k2) &
                       + dxdxi(1:3,3)*curlE(3,k2) ) / rjac

!        ...accumulate for the stiffness matrix: $((1/\mu) \tcurl E, \tcurl F)-((\omega^2 \eps - i \omega \sigma) E, F)$
            za = (CE(1)*CF(1) + CE(2)*CF(2) + CE(3)*CF(3)) / MU
            zb = (OMEGA*OMEGA*EPS - ZI*OMEGA*SIGMA) * (E(1)*F(1) + E(2)*F(2) + E(3)*F(3))
            Zaloc(k1,k2) = Zaloc(k1,k2) + (za-zb)*weight

   enddo; enddo; enddo
\end{lstlisting}

%--------------------------------------------------------------------
\subsection{DPG ultraweak implementation}
\label{sec:maxwell-ultraweak}

The broken ultraweak Maxwell formulation is given by:
\[
\left\{
\begin{split}
	\bs E, \bs H \in \bs{L^2}(\Omega), \hat{\bs E}, \hat{\bs H} \in H^{-1/2}(\tcurl, \Gammah): \bs n \times \bs E &= \bs n \times \bs E_0 \text{ on } \Gamma \, , \\
	(\bs H, \hcurl \bs F) - ((i \omega \eps + \sigma) \bs E , \bs F)
	&= (\bs J^{\text{imp}}, \bs F) &\quad \bs F \in \hHcurl \, , \\
	(\bs E, \hcurl \bs G) + \lb \bs n \times \hat{\bs E}, \bs G \rb_{\Gammah} + (i \omega \mu \bs H, \bs G)
	&= \bs 0 &\quad \bs G \in \hHcurl \, .
\end{split}
\right .
\]

Before looking at this section, the reader is encouraged to review the simpler DPG implementations of the Poisson problem, given in Sections~\ref{sec:poisson-primal} (primal) and \ref{sec:poisson-ultraweak} (ultraweak), which are discussed in greater detail. As in the Poisson DPG implementation, the static condensation of the Riesz representation of the residual is done at the element level (cf.~Section~\ref{sec:poisson-primal}). This model problem implementation uses a scaled adjoint graph norm for the ultraweak Maxwell formulation:
\[
	\| ( \bs F, \bs G ) \|^2_\test :=
	\| \hcurl \bs F - i \omega \bar{\mu} \bs G \|^2 +
	\| \hcurl \bs G + (i \omega \bar{\eps} - \sigma) \bs F \|^2 +
	\alpha ( \| \bs F \|^2 + \| \bs G \|^2) \, ,
\]
where $\bar{\cdot}$ indicates complex-conjugate, and $\alpha > 0$ is a scaling constant. See, e.g., \cite{melenk2023waveguide1, demkowicz2023waveguide2}, for both theoretical considerations and practical implications in terms of stability of choosing a proper value for $\alpha$.

The implementation is provided in \file{problems/MAXWELL/ULTRAWEAK\_DPG}.

\begin{lstlisting}[caption=\file{MAXWELL/ULTRAWEAK\_DPG/input/physics} input file.]
100000              MAXNODS, nodes anticipated
2                   NR_PHYSA, physics attributes
EHtrc   tangen 2    traces of electric/magnetic fields, H(curl), 2 components
EHfld   discon 6    electric/magnetic fields, L2, 3+3 components
\end{lstlisting}

In the ultraweak formulation, there are four physics unknowns: $\hat{\bs E}$, $\hat{\bs H}$, $\bs E$, $\bs H$; however, since all physics variables are solved in the same formulation and the electric and magnetic fields (and their traces) each use the same respective energy spaces (namely $\bs{L^2}$ and $H^{-1/2}(\tcurl)$), one can alternatively group these variables by specifying only two physics variables with double the number of components each: $(\hat{\bs E},\hat {\bs H}) \in (H^{-1/2}(\tcurl, \Gammah))^2$ and $(\bs E, \bs H) \in (\bs{L^2}(\Omega))^2$---recall that $\bs{L^2}(\Omega) = (L^2(\Omega))^3$. The variables are defined in the \file{physics} file and the tangential traces $(\hat{\bs E},\hat {\bs H})$ are specified as such by setting \code{\var{PHYSAi(1)}=.true.}. Choosing to specify grouped physics variables in such a way has certain implications for the ordering of degrees of freedom and the allocation of element-local matrices \var{ALOC} and \var{BLOC} (see Section~\ref{sec:coupled-variables}).

The element integration routine for the ultraweak DPG formulation is given by the following code:
\begin{lstlisting}[mathescape,caption=\file{MAXWELL/ULTRAWEAK\_DPG/}\routine{elem\_maxwell}: element integration]
!  ...use the enriched order to set the quadrature
      INTEGRATION = NORD_ADD
      call set_3D_int_DPG(ntype,norder,norient_face, nrint,xiloc,waloc)

!  ...loop over integration points
      do l=1,nrint

!     ...coordinates and weight of this integration point
         xi(1:3)=xiloc(1:3,l); wa=waloc(l)

!     ...H1 shape functions (for geometry)
         call shape3DH(ntype,xi,norder,norient_edge,norient_face, nrdofH,shapH,gradH)

!     ...L2 shape functions for the trial space
         call shape3DQ(ntype,xi,norder, nrdofQ,shapQ)

!     ...broken H(curl) shape functions for the enriched test space
         call shape3EE(ntype,xi,nordP, nrdofEE,shapEE,curlEE)

!     ...geometry map
         call geom3D(Mdle,xi,xnod,shapH,gradH,nrdofH, x,dxdxi,dxidx,rjac,iflag)

!     ...get permittivity at x
         call get_permittivity(mdle,x, eps)

!     ...integration weight
         weight = rjac*wa

!     ...get the RHS
         call getf(Mdle,x, zJ)

!     ...permittivity
         za = (ZI*OMEGA*EPS) * eps(:,:)

!     ...scalar permeability
         zc1 = ZI*OMEGA*MU

!     ...apply pullbacks
         call DGEMM('T','N',3,nrdofEE,3,1.d0     ,dxidx,3,shapEE,3,0.d0,shapF,3)
         call DGEMM('N','N',3,nrdofEE,3,1.d0/rjac,dxdxi,3,curlEE,3,0.d0,curlF,3)

!     ...apply permittivity
         zshapF = cmplx(shapF,0.d0,8)
         zcurlF = cmplx(curlF,0.d0,8)
         call ZGEMM('C','N',3,nrdofEE,3,ZONE,za,3,zshapF,3,ZERO,epsTshapF,3)
         call ZGEMM('N','N',3,nrdofEE,3,ZONE,za,3,zcurlF,3,ZERO,epscurlF ,3)

!     ...loop through enriched H(curl) test functions
         do k1=1,nrdofEE

!        ...pickup pulled-back test functions
            fldF(:) = shapF(:,k1);  crlF(:) = curlF(:,k1)
            fldG(:) = fldF(:);      crlG(:) = crlF(:)
            epsTfldF(:) = epsTshapF(:,k1)

!  --- load ---
!           $(J^{\text{imp}}, F)$ first equation RHS (with first H(curl) test function F)
!           $(0, G)$ second equation RHS is zero
            n = 2*k1-1
            bload_E(n) = bload_E(n) + (fldF(1)*zJ(1)+fldF(2)*zJ(2)+fldF(3)*zJ(3)) * weight

!  --- stiffness matrix ---
!        ...loop through L2 trial shape functions
            do k2=1,nrdofQ
!           ...first L2 variable
               m = (k2-1)*6
!           ...Piola transformation
               fldE(1:3) = shapQ(k2)/rjac; fldH = fldE

!           ...$-i \omega \eps (E,F)$
!           ...$(H,\curl F)$
               n = 2*k1-1
               stiff_EQ_T(m+1:m+3,n) = stiff_EQ_T(m+1:m+3,n) - fldE(:)*conjg(epsTfldF(:))*weight
               stiff_EQ_T(m+4:m+6,n) = stiff_EQ_T(m+4:m+6,n) + fldH(:)*crlF(:)*weight

!           ...$(E, \curl G)$
!           ...$i \omega \mu (H,G)$
               n = 2*k1
               stiff_EQ_T(m+1:m+3,n) = stiff_EQ_T(m+1:m+3,n) + fldE(:)*crlG(:)*weight
               stiff_EQ_T(m+4:m+6,n) = stiff_EQ_T(m+4:m+6,n) + zc1*fldH(:)*fldG(:)*weight
            enddo !..end of loop through L2 trial functions

!  --- Gram matrix ---
!        ...loop through enriched H(curl) test functions
            do k2=k1,nrdofEE
               fldE(:) = shapF(:,k2); epsTfldE(:) = epsTshapF(:,k2)
               crlE(:) = curlF(:,k2); epscrlE(:) = epscurlF(:,k2)

               call dot_product(fldF,fldE, FF)
               call dot_product(crlF,crlE, CC)

!          ...accumulate for the Hermitian Gram matrix (compute upper triangular only)
!             ------------------------
!             | (F_i,F_j)   (F_i,G_j) |    F_i/G_i are outer loop shape functions (fldF)
!             | (G_i,F_j)   (G_i,G_j) |    F_j/G_j are inner loop shape functions (fldE)
!             ------------------------

!             (F_j,F_i) terms = Int[F_^*i F_j] terms (G_11)
               n = 2*k1-1; m = 2*k2-1; k = nk(n,m)
               zaux = conjg(epsTfldF(1))*epsTfldE(1) + &
                       conjg(epsTfldF(2))*epsTfldE(2) + &
                       conjg(epsTfldF(3))*epsTfldE(3)
               gramP(k) = gramP(k) + (zaux + ALPHA_NORM*FF + CC)*weight

!              (G_j,F_i) terms = Int[F_^*i G_j] terms (G_12)
               n = 2*k1-1; m = 2*k2; k = nk(n,m)
               zaux = -(fldF(1)*epscrlE(1) + fldF(2)*epscrlE(2) + fldF(3)*epscrlE(3))
               zcux = conjg(zc1)*(crlF(1)*fldE(1) + crlF(2)*fldE(2) + crlF(3)*fldE(3))
               gramP(k) = gramP(k) + (zaux+zcux)*weight

!           ...compute lower triangular part of 2x2 G_ij matrix
!              only if it is not a diagonal element, G_ii
               if (k1 .ne. k2) then
!                 (F_j,G_i) terms = Int[G_^*i F_j] terms (G_21)
                  n = 2*k1; m = 2*k2-1; k = nk(n,m)
                  zaux = -(crlF(1)*epsTfldE(1) + crlF(2)*epsTfldE(2) + crlF(3)*epsTfldE(3))
                  zcux = zc1*(fldF(1)*crlE(1) + fldF(2)*crlE(2) + fldF(3)*crlE(3) )
                  gramP(k) = gramP(k) + (zaux+zcux)*weight
               endif

!              (G_j,G_i) terms = Int[G_^*i G_j] terms (G_22)
               n = 2*k1; m = 2*k2; k = nk(n,m)
               zcux = abs(zc1)**2*(fldF(1)*fldE(1) + fldF(2)*fldE(2) + fldF(3)*fldE(3))
               gramP(k) = gramP(k) + (zcux + ALPHA_NORM*FF + CC)*weight
         enddo; enddo !..end of loop through enriched H(curl) test functions
      enddo !..end of loop through integration points
\end{lstlisting}

Next, we provide the code performing the boundary integration:
\begin{lstlisting}[mathescape,caption=\file{MAXWELL/ULTRAWEAK\_DPG/}\routine{elem\_maxwell}: boundary integration]
!  ...loop through element faces
      do ifc=1,nrf

!     ...sign factor to determine the outward normal unit vector
         nsign = nsign_param(ntype,ifc)

!     ...face type
         ftype = face_type(ntype,ifc)

!     ...face order of approximation
         call face_order(ntype,ifc,norder, norderf)

!     ...set 2D quadrature
         INTEGRATION = NORD_ADD
         call set_2D_int_DPG(ftype,norderf,norient_face(ifc), nrint,tloc,wtloc)

!     ...loop through integration points
         do l=1,nrint

!        ...face coordinates
            t(1:2) = tloc(1:2,l)

!        ...face parametrization
            call face_param(ntype,ifc,t, xi,dxidt)

!        ...determine discontinuous H(curl) shape functions
            call shape3EE(ntype,xi,nordP, nrdof,shapEE,curlEE)

!        ...determine element H1 shape functions (for geometry)
            call shape3DH(ntype,xi,norder,norient_edge,norient_face, nrdof,shapH,gradH)

!        ...determine element H(curl) shape functions (for fluxes) on face
            call shape3DE(ntype,xi,norderi,norient_edge,norient_face, nrdof,shapE,curlE)

!        ...geometry map
            call bgeom3D(Mdle,xi,xnod,shapH,gradH,NrdofH,dxidt,nsign, &
                          x,dxdxi,dxidx,rjac,dxdt,rn,bjac)
            weight = bjac*wtloc(l)

!        ...pullback trial and test functions
            call DGEMM('T','N',3,nrdofEE,3,1.d0,dxidx,3,shapEE,3,0.d0,shapF ,3)
            call DGEMM('T','N',3,nrdofEi,3,1.d0,dxidx,3,shapE ,3,0.d0,shapFi,3)

!        ...loop through enriched H(curl) test functions
            do k1=1,nrdofEE
               E1(1:3) = shapF(:,k1)

!           ...loop through H(curl) trial functions
               do k2=1,NrdofEi
                  E2(1:3) = shapFi(:,k2)
                  call cross_product(rn,E2, rntimesE)

                  stiff_EE_T(2*k2-1,2*k1) = stiff_EE_T(2*k2-1,2*k1) + &
                              weight*( E1(1)*rntimesE(1) + E1(2)*rntimesE(2) + E1(3)*rntimesE(3) )
               enddo !..end loop through H(curl) trial functions
            enddo !..end loop through the enriched H(curl) test functions
         enddo !..end loop through integration points
      enddo !..end loop through element faces
\end{lstlisting}

Finally, the code performing the construction of the DPG linear system:
\begin{lstlisting}[mathescape,caption=\file{MAXWELL/ULTRAWEAK\_DPG/}\routine{elem\_maxwell}: constructing DPG linear system.]
      allocate(stiff_ALL(NrTest,NrTrial+1))

!  ...Total test/trial DOFs of the element
      i1 = NrTest ; j1 = 2*NrdofEi ; j2 = 6*NrdofQ

! ...Copy stiffness and load into one matrix
      stiff_ALL(1:i1,1:j1) = transpose(stiff_EE_T(1:j1,1:i1))
      stiff_ALL(1:i1,j1+1:j1+j2) = transpose(stiff_EQ_T(1:j2,1:i1))
      stiff_ALL(1:i1,j1+j2+1) = bload_E(1:i1)

!  ...A. Compute Cholesky factorization of Gram Matrix, $\mr{G=U^* U (=LL^*)}$
      call ZPPTRF('U',NrTest,gramP,info)

!  ...B. Solve triangular system to obtain $\mr{\tilde{B}}$, $\mr{(LX=) U^* X = [B|l]}$
      call ZTPTRS('U','C','N',NrTest,NrTrial+1,gramP,stiff_ALL,NrTest,info)

      allocate(zaloc(NrTrial+1,NrTrial+1)); zaloc = ZERO

!  ...C. Matrix multiply: $\mr{B^* G^{-1} B (=\tilde{B}^* \tilde{B})}$
      call ZHERK('U','C',NrTrial+1,NrTest,ZONE,stiff_ALL,NrTest,ZERO,zaloc,NrTrial+1)

!  ...D. Fill lower triangular part of Hermitian matrix $\mr{\tilde B^* \tilde B}$
      do i=1,NrTrial
         zaloc(i+1:NrTrial+1,i) = conjg(zaloc(i,i+1:NrTrial+1))
      enddo

!..zaloc has now all blocks of the stiffness and load:
!  $z_{\text{aloc}} =  \var{ALOC(1,1)} \quad \var{ALOC(1,2)} \quad \var{BLOC(1)}$
!         $\, \var{ALOC(2,1)} \quad \var{ALOC(2,2)} \quad \var{BLOC(2)}$
\end{lstlisting}

%% file: Chapters/APPENDIX_B_APPLICATIONS/appendix_applications_chapter.tex
%
%!TEX root = ../../hp3D_user_guide.tex
%

\chapter{Applications}
\label{chap:applications}

%--------------------------------------------------------------------

To provide further examples and give the reader an idea of the scope of computations \hp3D is suitable for, this chapter briefly summarizes applications that have been implemented within the current version of the \hp3D finite element code. If an implementation is available in the public repository, we provide a file path for the interested user. References for further reading are provided in all cases.

%--------------------------------------------------------------------
\section{Optical Fiber Amplifier}
\label{sec:laser}

Fiber amplifiers are optical waveguides made of silica glass designed for power-scaling highly-coherent laser light. At high optical intensities, undesired nonlinear effects may negatively affect the beam quality. \hp3D has been used to study such nonlinear effects (e.g.~the interplay between the propagating electromagnetic fields and thermal effects) for active gain fiber amplifiers based on a finite element model of the time-harmonic Maxwell equations coupled with the heat equation \cite{henneking2021fiber,nagaraj2018raman}. An implementation of this application is available in the \file{problems/LASER} directory. 

In addition to the coupled multiphysics formulation, the application employs a high-order discretization and anisotropic adaptive refinements \cite{henneking2021pollution} for a hybrid mesh consisting of both prismatic and hexahedral elements with curvilinear geometry. The application also served as the first real testbed problem for large-scale computation with the MPI/OpenMP parallel \hp3D code, successfully scaling up to thousands of wavelengths and $\sim$1B degrees of freedom \cite{henneking2021phd,henneking2022parallel}. The DPG FE implementation of the Maxwell equations also features fast integration routines via sum factorization for both hexahedral and prismatic elements \cite{mora2019fast,badger2020fast}.

%--------------------------------------------------------------------
\section{Adaptive Solution of High-Frequency Acoustic and Electromagnetic Scattering}
\label{sec:adaptive}

Accurate and efficient numerical simulations of wave propagation phenomena are very crucial in numerous engineering and physics applications, such as non-destructive testing, plasma fusion, modeling of meta-materials and biomedical and radar imaging. With traditional discretization methods, these simulations are extremely challenging due to two major issues: instability and indefiniteness of the linear systems. Consequently, common cutting-edge elliptic solvers simply break down. The DPG method overcomes both issues. As a non-standard minimum residual method, it promises high accuracy, unconditional stability and definite linear systems.

This work implemented a novel multigrid solver within \hp3D for linear systems arising from DPG discretizations with a special focus in acoustic and electromagnetic wave propagation problems. The construction is heavily based on the attractive properties of the DPG method, but also on well-established theory of Schwarz domain decomposition and multigrid methods. As it is showcased in \cite{petrides2019phd,petrides2021adaptive,petrides2017adaptive}, the method is stable and reliable, and it is suitable for adaptive $hp$-meshes. The solver works hand-in-glove with the built-in DPG error indicator to drive adaptive mesh refinements. Integrating the iterative solver with the adaptive refinement procedure enables efficient and accurate solutions to challenging problems in the high-frequency regime. A distributed parallel version of the DPG multigrid solver is currently under development and is discussed in \cite{badger2023scalable}.

%--------------------------------------------------------------------
\section{Linear Elasticity with Two Materials}
\label{sec:hose}

Linear elasticity is a continuum model of the deformation and stresses of solids. The first finite element methods were developed to predict the effects of loads on structures modeled with linear elasticity and the study of linearly elastic materials has continued to be an important application of finite element analysis to this day. A finite element method for linear elasticity is available in the \file{problems/SHEATHED\_HOSE} directory.

The implementation employs distinct DPG methods (primal and ultraweak \cite{keith2016elasticity}) for two separate subdomains in curvilinear body. The primal DPG method is used in a subdomain with compressible material (steel) and the ultraweak DPG method is used in a subdomain with incompressible material (rubber) \cite{fuentes2017coupled}. The coupling between the two methods is naturally incorporated through the common trace and flux dofs on the element boundary. This application demonstrates support in hp3D for subdomain-dependent local element assembly, which can be used to support the solution of multiphysics problems.

%--------------------------------------------------------------------
\section{Nonlinear Elasticity}
\label{sec:nonlinear-elasticity}

Beyond linear elasticity problems, it is common in computational mechanics to model and simulate nonlinear elasticity cases, e.g.~to compute large deformations. In such cases, one needs to account for a non-convex stored energy functional that must be minimized. A practical approach for solving this minimization problem is to apply loads in incremental steps to find successive local minima via Newton--Raphson methods. Each iteration of the Newton--Raphson method applied to the system of nonlinear equations solves a linear system that is derived from a variational formulation of the linearized boundary value problem (see \cite{mora2020polydpg}, Chapters 6 and 7). 

The implementation of the nonlinear elasticity problem in \hp3D assumes an isotropic hyperelastic material with dead load undergoing quasi-static deformation. Various variational formulations of the linearized equations can be employed for this scenario and numerically solved in \hp3D with user-defined convergence tolerances. The application directory, \file{problems/NONLINEAR\_ELASTICITY}, includes several implementations: the Bubnov--Galerkin discretization of the classical variational formulation in \file{GALERKIN}; the DPG discretization of the broken primal formulation in \file{DPG\_PRIMAL}; and the DPG discretization of the broken ultraweak formulation in \file{DPG\_UW}. The mathematical theory and material constitutive models are detailed in \cite{mora2020polydpg}. A more thorough description of the implementation for this problem is presented in  \cite{hpbook3}.

%--------------------------------------------------------------------
\section{Time-Harmonic Applications in Linear Viscoelasticity}
\label{sec:visco-elasticity}

The linearized equations describing the dynamics of viscoelastic material behavior in the time-harmonic regime are very similar to those of time-harmonic linear elasticity. The difference is that the material properties themselves (as well as the relevant unknowns) are complex-valued. The reason is that the (linearized) constitutive model describing the relationship between stress and strain depends on the strain history, and is given by a Volterra-type convolution, which becomes a product in the time-harmonic regime. This allows to model real-life experiments in dynamic mechanical analysis (DMA) as well as phenomena that involve low-amplitude vibrations of viscoelastic material at different frequencies.

An example where \hp3D was used includes the faithful reproduction of a DMA experiment involving a thin beam. Here support for complex variables is needed, as well as anisotropic $p$ and local anisotropic refinements in $h$ when the domain is discretized with hexahedra. The reason is that adaptivity is localized in certain regions where most of the deformation gradient is concentrated due to nontrivial clamping effects \cite{fuentes2017viscoelasticity}. Moreover, in an application involving form-wound stator coils found in electric machinery, it is necessary to introduce a nontrivial ``surface'' forcing coming from high-frequency Lorentz forces. This type of forcing is not typically supported. However, with the proper variational formulation and boundary conditions, it is possible to use \hp3D to solve this problem, which also includes some of the features mentioned previously \cite{fuentes2018phd}.

%--------------------------------------------------------------------
\section{Acoustics of the Human Head}
\label{sec:human-head}

The problem was formulated and solved in the course of the dissertation of Paolo Gatto devoted to the analysis of transmission of sound in the human head \cite{gatto2012phd}, see also \cite{demkowicz2011bone,gatto2013bone}. More specifically, it studied the mechanism of transmitting exterior acoustic energy into the cochlea with the ear canal {\em blocked}. The domain of interest included the human skull modeled with a spherical shell with an inner structure representing the ear canal, a cochlea with the tympanic membrane, and a simplified model of ossicles consisting of a single elastic rod connecting the ear drum with the oval window in the cochlea. The skull was surrounded with an additional layer representing air, terminated with a PML.
The skull, ear canal, cochlea, ossicles and the three membranes were modeled with elasticity; whereas the air, exterior and interior ear, the brain and the fluid inside of the cochlea were modeled with acoustics, a total of 16 subdomains whose geometry was represented within the GMP package. The initial tetrahedral linear geometry/mesh was obtained using NETGEN and imported into GMP where two important modifications were made. The skull and all membranes, represented initially with zero thickness interfaces, were extruded into thin-walled structures, and the triangular mesh representing the skull was extruded into two additional layers for air and PML. The extrusions brought in both prismatic and pyramidal elements which motivated adding these elements to the \hp3D code.

In the elastic subdomains, the code supported the elastic displacement field---a single $H^1$ field with three components. In the acoustical subdomains, the code supported the pressure---a single-component $H^1$ field. Note that the nodes located on acoustic/elastic interfaces supported both variables. More precisely, in the external acoustic domain the acoustic variable represented the scattered pressure, and in the internal acoustical domains the total pressure. The forcing came from a plane wave impinging on the skull, resulting with an additional term on the exterior acoustics/elasticity interface. A standard variational formulation (see e.g.~\cite{hpbook2}) was used with couplings between the elasticity and acoustic problems enforced weakly by incorporating interface integrals in the variational formulation. 
%It is worth mentioning that the (asymptotic) discrete stability of the problem is guaranteed by the Mikhlin theory \cite{demkowicz2020fem}.
The problem was discretized with the standard Galerkin method. In the initial mesh, the three membranes were discretized with quintic elements to avoid locking, with the rest of the domain covered with quadratic elements except for the PML subdomain where higher order anisotropic elements were used. Original shape functions for $H^1$-conforming elements of all shapes were developed in \cite{gatto2010shape}. Finally, the problem was solved using $h$-adaptivity driven by a problem-dependent, implicit a-posteriori error estimate developed in the course of the project.

%--------------------------------------------------------------------
\section{Electromagnetic Radiation and Induced Heat Transfer in the Human Body}
\label{sec:human-body}

The project focused on estimating heating effects in the human head from the absorption of electromagnetic (EM) waves emanating from a  cell phone; it constituted (partially) the dissertation of Kyungjoo Kim \cite{Kim2013phd}.

The transient heat equation was coupled with time-harmonic Maxwell equations; the standard Galerkin method was used for both heat and Maxwell problem. The heat equation was discretized with the implicit Euler method. At each time step, the new temperature field was used to compute new (temperature-dependent) material constants for the Maxwell equations. Due to disparate time scales, the time-harmonic version of the Maxwell equations was employed. The resulting new electric field was then used to update the source in the heat equation. The iterations were continued until a steady-state was reached.

The head was modeled with both simplified multilayer spherical geometry as well as a more precise (homogeneous) head model imported from COMSOL. The application was the second problem driving the original development of the current version of the \hp3D code. The Maxwell problem was discretized with $H(\tcurl)$-conforming tetrahedral and prismatic elements.

\clearpage

%% file: Chapters/APPENDIX_C_DPG/appendix_dpg_chapter.tex
%
%!TEX root = ../../hp3D_user_guide.tex
%

\chapter{Implementation of the DPG Method}
\label{chap:dpg}

%--------------------------------------------------------------------

The \hp3D library has been the main tool for the numerical studies with the discontinuous Petrov--Galerkin (DPG) Method with Optimal Test Functions \cite{demkowicz2017dpg}. Many of the recent applications implemented in \hp3D (see Appendix~\ref{chap:applications}) were discretized with the DPG method. DPG implementations are straightforward to implement within the \hp3D framework, and we briefly review (in abstract form) how the DPG linear system is constructed. 

Examples of model problem implementations with DPG are also provided (e.g.~Appendix~\ref{sec:poisson-primal}).

\paragraph{The DPG mixed formulation.} Consider an abstract variational problem and its Petrov--Galerkin (PG) discretization:
\[
\left\{
\begin{array}{lllll}
u \in \trial \\
b(u,v) = l(v) \quad v \in \test
\end{array}
\right.
\qquad
\longrightarrow
\qquad
\left\{
\begin{array}{lllll}
u_h \in \trial_h \subset \trial \\
b(u_h,v_h) = l(v_h) \quad v_h \in \test_h \subset \test
\end{array}
\right.
\]
The {\em Petrov--Galerkin } discretization of the problem requires $\dim \trial_h = \dim \test_h$. For $\trial=\test$, the choice of $\trial_h=\test_h$ leads to the {\em Bubnov--Galerkin} (BG) method.
In the {\em Petrov--Galerkin Method with Optimal Test Functions}, we embed the original problem into an equivalent {\em mixed problem}:
\[
\arraycolsep=2pt
\left\{
\begin{array}{lllll}
\psi \in \test, u \in \trial \\
(\psi,v)_\test + b(u,v) & = l(v)  \quad &\ v \in \test\\
b(w,\psi) & = 0 &\ w \in \trial
\end{array}
\right.
\qquad
\longrightarrow
\qquad
\left\{
\begin{array}{lllll}
\psi_h \in \test_h ,\, u_h \in \trial_h \\
(\psi_h,v_h)_\test + b(u_h,v_h) & = l(v_h) \quad &\ v_h \in \test_h \\
b(w_h, \psi_h) & = 0 &\ w_h \in \trial_h
\end{array}
\right.
\]
The mixed formulation involves an additional unknown---$\psi \in \test$---the Riesz representation of the residual. On the continuous level $\psi = 0$ but its discretization $\psi_h \in \test_h$ is different from zero. The test inner product $(\cdot,\cdot)_\test$ enters directly the problem and affects the corresponding FE solution $u_h$. The main motivation for solving the more expensive mixed problem comes from stability considerations: the discrete spaces $\trial_h,\test_h$ need not longer be of the same dimension, and it is easier to satisfy the \emph{discrete inf--sup condition} by simply employing a larger discrete test space. We are solving effectively an overdetermined discrete system.

An alternative is to test with a larger space of {\em broken} test functions $V(\Omega_h) \supset \test$ where $\Omega_h$ denotes the 
decomposition of the domain into finite elements. Testing from a larger space of broken test functions necessitates the introduction of an 
additional unknown---the {\em trace} $\hat{u}$---a Lagrange multiplier. The modified problem looks as follows:
\[
\left\{
\begin{array}{llllll}
u \in \trial ,\, \hat{u} \in \hat{\trial} \\
\underbrace{b(u,v) + \langle \hat{u}, v \rangle_{\Gamma_h}}_{b_{\text{mod}}((u,\hat{u}),v)} = l(v) \quad v \in \test(\Omega_h)\, .
\end{array}
\right.
\]
The new group variable includes the original unknown and the trace unknown. The corresponding PG scheme with optimal test functions is known as the DPG method. The word {\em discontinuous} in the name refers to the use of broken (discontinuous) test spaces.
\[
\arraycolsep=2pt
\left\{
\begin{array}{llllll}
u \in \trial ,\, \hat{u} \in \hat{\trial}, \psi \in \test(\Omega_h) \\
(\psi,v)_{\test(\Omega_h)} + b(u,v) + \langle \hat{u}, v \rangle_{\Gamma_h} & = l(v) \quad &\quad v \in \test(\Omega_h)\\
b(w, \psi) & = 0 &\quad w \in \trial\\
\lb \hat{w}, v \rb_{\Gamma_h} & = 0 &\quad \hat{w} \in \hat{\trial} \, .
\end{array}
\right.
\]

The main advantage of using the broken test space in context of the PG method with optimal test functions is that the Gram matrix resulting from the discretization of
the test inner product $(\psi,v)_\test$ becomes block-diagonal, and the additional variable $\psi \in V(\Omega_h)$ can be statically condensed on the
element level. The whole ``DPG magic'' happens only in the element routine, the remaining part of the code remains the same as for the standard Bubnov--Galerkin method. The mixed problem,
\[
\left(
\begin{array}{ccc}
G & B & \hat{B} \\
B^* & 0& 0\\
\hat{B}^* & 0& 0
\end{array}
\right)
\, 
\left(
\begin{array}{c}
\psi\\u\\\hat{u}
\end{array}
\right)
\, = \,
\left(
\begin{array}{c}
l \\0 \\0
\end{array}
\right)\, ,
\]
becomes a standard Bubnov--Galerkin problem:
\[
\left(
\begin{array}{cc}
B^* G^{-1} B & B^* G^{-1} \hat{B} \\
\hat{B}^* G^{-1} B^* & \hat{B}^* G^{-1} \hat{B}^*
\end{array}
\right)
\, 
\left(
\begin{array}{c}
u\\\hat{u}
\end{array}
\right)
\, = \,
\left(
\begin{array}{c}
B^* G^{-1}  l \\\hat{B}^* G^{-1} l
\end{array}
\right)\, .
\]

Recall that the trace unknown $\hat u$ is discretized in \hp3D by using element-wise restrictions of $H^1$-, $H(\tcurl)$-, and $H(\tdiv)$-conforming elements. This is accomplished by first defining a standard $H^1$, $H(\tcurl)$, or $H(\tdiv)$ physics variable in the \file{physics} input file and then declaring it a trace via the global flag array \var{PHYSAi(:)}. For more information on the traces, see Section~\ref{sec:traces}.

\paragraph{Constructing the DPG linear system.}
In practice, the local element matrices of the DPG linear system can be constructed in the following way:

First, reduce the broken mixed formulation to a matrix equation:
\begin{itemize}
\item{
Let $\mf{U}_h = \{ \mf{u}_i \}_{i=1}^N,\ \mf{\hat U}_h = \{ \mf{\hat u}_i \}_{i=1}^{\hat N}$, and $\mf{V}^r = \{ \mf{v}_i \}_{i=1}^M$ (where $M > N+\hat N$) denote bases for the discrete trial space $\trial_h \times \hat \trial_h$ and the enriched test space $V^r$, respectively.}
\item{
Define the stiffness matrix for the modified bilinear form, the Gram matrix, and the load vector:
\[
	\mr{B}_{ij} = b(\mf{u}_j, \mf{v}_i), \quad 
	\mr{\hat B}_{ij} = \lb \mf{\hat u}_j, \mf{v}_i \rb, \quad 
	\mr{G}_{ij} = (\mf{v}_j, \mf{v}_i)_V, \quad 
	\mr l_i = l(\mf{v}_i) .
\]
}
\item{
In matrix form, the problem can now be formulated in the following way:

Find the set of coefficients (over field $\bb{F} = \bb{R}$ (or $\bb{C}$))
\[
	\mr w = [\mr{w}_i]_{i=1}^N \in \bb{F}^N, \quad
	\mr{\hat w} \in [\mr{\hat w}_i]_{i=1}^{\hat N} \in \bb{F}^{\hat N}, \text{ and }\
	\mr{q} = [\mr{q}_i]_{i=1}^M \in \bb{F}^M
\]
such that
\[
	\mr {u_h} = \sum_{i=1}^{N} \mr{w}_i \mf{u}_i , \quad
	\mr{\hat{u}_h} = \sum_{i=1}^{\hat N} \mr{\hat w}_i \mf{\hat u}_i , \text{ and }\
	\mr{\Psi} = \sum_{i=1}^{M} \mr{q}_i \mf{v}_i
\]
satisfy
\[
	\left[ \begin{array}{ccc}
		\mr G & \ \mr B \ \ & \mr{\hat B} \\
		\mr B^* & \ \mr 0 \ \ & \mr 0\\
		\mr{\hat B^*} & \ \mr 0 \ \ & \mr 0
	\end{array} \right]
	\left[ \begin{array}{c}
		\mr{\Psi} \\
		\mr{u_h} \\
		\mr{\hat u_h}
	\end{array} \right]
	=
	\left[ \begin{array}{c}
		\mr l \\
		\mr 0 \\
		\mr 0
	\end{array} \right] .
\]
}
\end{itemize}

Let $\mr{B} := [\mr{B}\ \mr{\hat B}]$. Then, \vskip -30pt
\begin{align*}
	\mr{B^*} \mr{G^{-1}} \mr{B} &= 
	\mr{B^*} \mr{(L L^*)^{-1}} \mr{B} =
	( \mr{L^{-1}} \mr{B} )^* \overbrace{( \mr{L^{-1}} \mr{B} )}^{\mr{\tilde B}} =
	\mr{\tilde B}^* \mr{\tilde B} \, ,
	\\[3pt]
	\mr{B^*} \mr{G^{-1}} \mr{l} &=
	( \mr{L^{-1}} \mr{B} )^* ( \mr{L^{-1} l} ) =
	\mr{\tilde B}^* ( \mr{L^{-1} l} ) \, .
\end{align*}
That is, the statically condensed DPG linear system is computed as follows:
\begin{enumerate}
	\itemsep 0pt
	\item Cholesky factorization: \hskip 7pt 
	$\mr{G = L L^*}$
	\item Solve triangular system: \hskip 2pt 
	$\mr{[\tilde B \, | \, \tilde l \,] = L^{-1} [B \, | \, l \,]}$
	\item Matrix multiply: \hskip 37pt 
	$\mr{\tilde B^* [\tilde B \, | \, \tilde l \,]}$
\end{enumerate}
The \routine{elem} routine concludes by storing this statically condensed system into the corresponding blocks of the assembly arrays \var{ALOC} and \var{BLOC} (cf.~\ref{sec:coupled-variables}).

%\input{Chapters/APPENDIX_B_APPLICATIONS/comments}